\documentclass[12pt,a4paper,oneside,notitlepage]{amsart}
\usepackage{amssymb}
\usepackage[arrow,matrix,curve]{xy}\SilentMatrices
\def\xyma{\xymatrix@M.7em}
\usepackage[centertags]{amsmath}
\usepackage[T1]{fontenc}
\usepackage[latin1]{inputenc}

\usepackage{amsfonts}
\usepackage{amscd}
\usepackage{euscript}
\usepackage{amsmath}

\numberwithin{equation}{section}
\newtheorem{theorem}{Theorem}[section]

\newtheorem{defi}[theorem]{Definition}
\newtheorem{prop}[theorem]{Proposition}
\newtheorem{lemma}[theorem]{Lemma}
\newtheorem{cor}[theorem]{Corollary}
\newcommand{\Para}{\par\vspace{.5cm}\noindent}
\newcommand{\para}{\par\vspace{.5cm}}
\newcommand{\Z}{\mathbb Z}
\newcommand{\Q}{\mathbb Q}
\usepackage{amsfonts,amssymb}

\newtheorem{satz}{Theorem}[section]

\newtheorem{kor}[satz]{Corollary}
\newtheorem{lem}[satz]{Lemma}

\newtheorem{bsp}[satz]{Example}

\newcounter{Roma}
\newcounter{Ara}
\newcounter{let}

\makeindex

\newcommand{\nc}{\newcommand}

\nc{\mapco}{\,\colon\, }
\nc{\ab}{^{ab}}

\nc{\catc}{{\C}}
\nc{\we}{\vee}

\nc{\hra}{\hookrightarrow}

\nc{\epi}{epimorphism}
\nc{\repi}{regular epimorphism}
\nc{\mono}{monomorphism}
\nc{\iso}{isomorphism}
\nc{\coker}[1]{\mbox{${\rm Coker}(#1)$}}
\nc{\Ker}[1]{\mbox{Ker$(#1)$}}

\nc{\defgl}{\stackrel{def}{=}}

\nc{\V}{\vspace{5mm}}
\nc{\VV}{\vspace{4mm}}
\nc{\lra}{\longrightarrow}
\nc{\lla}{\longleftarrow}
\nc{\mr}[1]{ \stackrel{#1}{\lra} }
\nc{\ml}[1]{ \stackrel{#1}{\lla} }
\nc{\hmr}[1]{\hspace{2mm}\stackrel{#1}{\lra}\hspace{2mm}}
\nc{\hml}[1]{\hspace{2mm} \stackrel{#1}{\lla}\hspace{2mm}}
\nc{\N}{\noindent}
\nc{\st}{^{\prime}}
\nc{\ot}{\otimes}
\nc{\hcong}{ \hspace{2mm}\cong\hspace{2mm}  }
\nc{\hfbox}{\hfill$\Box$}
\nc{\REF}[1]{(\ref{#1})}

\def\Z{\ifmmode{Z\hskip -4.8pt Z} \else{\hbox{$Z\hskip -4.8pt Z$}}\fi}

\def\Q{\ifmmode{Q\hskip-5.0pt\vrule height6.0pt depth
0pt\hskip6pt}\else{\hbox{$Q\hskip-5.0pt\vrule height6.0pt depth
0pt\hskip6pt$}}\fi}

\nc{\Ph}{\phantom{}}

\nc{\BE}{\begin{equation}}
\nc{\EE}{\end{equation}}

\nc{\dst}{\displaystyle}
\nc{\sst}{\scriptscriptstyle}
\nc{\ssst}{\scriptscriptstyle}
\nc{\proofof}[1]{\N{\bf Proof of {#1}\,:}\quad}
\nc{\proofofthm}[1]{\N{\bf Proof of theorem \ref{#1}\,:}\quad}

\nc{\htt}[1]{^{\otimes #1}}

\nc{\Sur}[1]{\mbox{$\:\stackrel{#1}{\lra\!\!\!\!\!\to\,}\:$}}

\def\INJ{\mbox{\mathsurround=0pt
\makebox[0mm][r]{\parbox{0mm}{\rule[-0.65mm]{0mm}{0.2mm}$\scriptscriptstyle>$}}
\makebox[0.7cm][l]{\parbox{0.7cm} {$\lra$}}}}

\def\Inj#1{\mbox{$\:\stackrel{#1}{\INJ}\:$}}

\nc{\Injup}[1]{\mapup{#1}}
\nc{\injup}[1]{\mapup{#1}}
\nc{\injdown}[1]{\mapdown{#1}}

\def\mapup#1{\mbox{$\rule[-1mm]{0cm}{0.7cm}  
\makebox[0mm][r]{\raisebox{0.2mm}{$\scriptstyle\phantom{\cong}$}\hspace{0.6mm}}
\bigg\uparrow\rlap{$\vcenter{\hbox{$\scriptstyle#1$}}$}$}}

\def\mapdown#1{\mbox{$\rule[-1mm]{0cm}{0.7cm}  
\makebox[0mm][r]{\raisebox{0.2mm}{$\scriptstyle $}\hspace{0.6mm}}
\bigg\downarrow\rlap{$\vcenter{\hbox{$\scriptstyle#1$}}$}$}}

\def\isor#1{\mbox{$\smash{\mathop{\longrightarrow}\limits^{\cong}_{#1}}$}}

\def\surdown#1{\makebox[0mm]{$\mapdown{\raisebox{0.4mm}{$\scriptstyle#1$}}$}
\makebox[0mm]{\raisebox{-2.15mm}{$\downarrow$}}}

\def\surup#1{\makebox[0mm]{$\mapup{\raisebox{0.4mm}{$\scriptstyle#1$}}$}
\makebox[0mm]{\raisebox{0.7mm}{$\mapup{}$}}}

\newcommand{\brokenr}[2]{
@\cdhgeneric>\raise2.8pt\hbox to3.5pt{\hrulefill}\mkern9mu>
\raise2.8pt\hbox to3.5pt{\hrulefill}\hbox to5pt{}>
\mkern-7.5mu\dashrightarrow>#1 >
#2>}

\newcommand{\brokenup}[1]{
@\cdvgeneric>\hat\cdot
>\raisebox{3pt}{$\vdots$}>
\vbox{\kern3.5pt\hbox{$\cdot$}\kern-3.5pt}
> >\hspace{1mm}#1>}

\newcommand{\functlr}[2]{
\raisebox{0.4pt}{$\hss\begin{CD}
@>\vbox{\hbox to 0pt{$\hss\begin{CD}@<#1<<\end{CD}\hss$}\vskip-2pt}
>#2 >
\end{CD}\hss$}
}

\newcommand{\maprr}[2]{
\raisebox{-0.9pt}{$\hss\begin{CD}
@>\vbox{\hbox to 0pt{$\hss\begin{CD}@>#1>>\end{CD}\hss$}\vskip-3pt}
>#2 >
\end{CD}\hss$}
}

\newcommand{\mapdd}[2]{
@\cdvstandard>\downarrow\hspace{1.5pt}\kern-4pt\hspace{1.5pt}\Big\downarrow>#1>#2>}

\newcommand{\mapud}[2]{
@\cdvstandard>\uparrow\kern-3pt\Big\downarrow>#1>#2>}

\nc{\auf}{\twoheadrightarrow}
 
\nc{\ruled}{\rule[-4mm]{0mm}{0mm}}

\nc{\qu}{quadratic}
\nc{\GBoGB}{G/BG\st \otimes G/BG\st}
\nc{\lstar}{_{\raisebox{-1mm}{$*$}}}
\nc{\sm}{\:{ \wedge}\:}

\nc{\Rmod}{${\bf R}$-module}

\nc{\map}[3]{\mbox{$#1 \mapco #2 \to #3$}}

\nc{\rond}{{\,\sst \circ\,}}

\nc{\ruleu}{\rule{0mm}{7mm}}

\nc{\T}[1]{\tilde{#1}}

\nc{\Imm}[1]{\mbox{${\rm Im}(#1)$}}

\nc{\tw}{\end{document}}

\nc{\QJG}{\frac{\dst I(G) J}{\dst I^2(G) J}}
\nc{\QJH}{\frac{\dst I(H) J}{\dst I^2(H) J}}
\nc{\otz}{\ot}
\nc{\UL}[2]{{\rm U}_{#1}{\rm L}(#2)}

\nc{\IRN}[1]{I_{R,\cal G}^{#1}(G)}

\nc{\ULG}[1]{{\rm U}_{#1}{\rm L}^{\cal G}(G)}
\nc{\UGH}[1]{{\rm U}_{#1}^{\cal G}(G,H)}

\nc{\LG}[1]{{\rm L}^{\cal G}_{#1}(G)}

\nc{\calG}{{\cal G}}
\nc{\AB}{^{AB}}

 \nc{\IE}[1]{\mbox{$I^{#1}(E)$}}
 \nc{\IN}[1]{\mbox{$I^{#1}(N)$}}
\nc{\IG}[1]{\mbox{$I^{#1}(G)$}}

\begin{document}
\begin{center}{\bf \Large Dimension Quotients}\end{center}
\par\vspace{2cm}\centerline{M. Hartl, R. Mikhailov\footnote{The research of the second  author was supported by Russian Presidential Grant No.
MK- 3466.2007.1, and Russian Science
Support Foundation} and I. B. S.
Passi\footnote{INSA Senior Scientist}}
\par\vspace{2cm}\begin{small}\begin{quote}\centerline{\sf
Abstract}\para We present two approaches, one homological and the
other simplicial,  for the investigation of dimension quotients
of groups.
 The theory is illustrated, in particular,  with a conceptual discussion of the fourth and fifth dimension quotients.
\end{quote}
\end{small}
\para
\begin{tiny}\begin{quote}\tableofcontents\end{quote}\end{tiny}
\Para
\section{Introduction}

\para
Let $G$ be a group, $\mathbb Z[G]$ its integral group ring and $\mathfrak g$ the augmentation ideal. The dimension subgroups $D_n(G),\ n\geq 1,$ of $G$ are the subgroups determined by the augmentation powers $\mathfrak g^n$, i.e., $D_n(G)=G\cap (1+\mathfrak g^n)$. The dimension subgroups were first defined by W. Magnus \cite{Magnus:37} for free groups, via power series ring in non-commuting variables,   and ever since they have been a subject of intensive study (see \cite{Passi:79}, \cite{Gupta:87}). It is easy to see that, for all $n\geq 1$, $D_n(G)$ contains $\gamma_n(G)$, the $n$-th term in the lower central series of $G$ and equality holds for $n=1,\,2,\,3$. 
Refuting what came to be known as the {\it dimension conjecture}, E. Rips \cite{Rips:72} constructed a nilpotent 2-group of class three having  non-trivial fourth dimension subgroup.  In a significant contribution J. A. Sj\"{o}gren \cite{Sjogren:79} proved that there exist constants $c_n,\ n\geq 1$, such that, for every group $G$, the exponent of the quotient $D_n(G)/\gamma_n(G)$, called the $n$-th dimension quotient of $G$, divides $c_n$.   While considerable work has been done on the dimension  subgroups, notably by N. Gupta (\cite{Gupta:90a}, \cite{Gupta:91b}, \cite{Gupta:02}), the structure of the dimension quotients, in general, remains intractable. 
\para
The purpose of this paper is to present two approaches to study the dimension quotients, one homological and the other using simplicial objects. The homological approach originates with \cite{Passi:68b}, where it is shown that for odd prime power groups the fourth dimension quotient is trivial. The simplicial approach rests on the definition of derived functors of arbitrary endo-functors on the category of groups and certain spectral sequences studied by E. B. Curtis (\cite{Curtis:63}, \cite{Curtis:71}) in the context of simplicial homotopy theory. We illustrate our methods with a conceptual discussion of the fourth and fifth dimension quotients. We hope that the methods outlined here will open fresh avenues to investigate the dimension quotients. The main point to be observed in our approach is that it brings out new relationship between dimension subgroups, derived functors of non-additive functors and simplicial homotopy theory, thus providing new tools and connections. \para 
We begin by recalling in section \ref{polfunct}
the constructions of certain quadratic functors that are needed in this work. We then take up in section \ref{relative1} the study  of relative dimension subgroups and polynomial (co)homology. For any normal subgroup $N$ of a group $E$, the $n$-th dimension subgroup $ D_n(E,\,N)$ of $E$ relative to $N$ is defined by 
\[ D_n(E,\,N) := E \cap \Big( 1 + \mathfrak n\mathfrak e + \mathfrak e^n \Big)\;.\]
Let $R$ be a normal subgroup of a free group $F$. We prove (Theorem \ref{D3free}) that there exists a   natural isomorphism
\[   (F/R\gamma_2(F))\stackrel{\sst\sm}{*} (F/R\gamma_2(F))\hcong D_3(F,\,R)/\gamma_2(R)\gamma_3(F) \,,\] where, for an Abelian group $A$, $A\stackrel{\sst\sm}{*}A$ denotes the exterior torsion square of $A$ [see (\ref{LSdef}) for definition]. An application of our Theorem \ref{D3free} to group cohomology $H^2_{{\mathcal N}_c}(G,\,M)$ with respect to the variety ${\mathcal N}_c$
of $c$-step nilpotent groups, for $c=2$, leads to a solution (Theorem \ref{Kuc=2}) of  a problem of Leedham--Green \cite{LG:71} in this case. 
Let $e\mapco  N\stackrel{i}{\hookrightarrow} E \stackrel{q}{\auf} G$ be a central extension of groups. We give a homological characterisation of $D_n(E,N)\cap N$ (see Theorem \ref{DnH2}). In section \ref{fourthdq1},
restricting our analysis to the case $n=4$ (see Theorem \ref{Kerrho2}), we go on to provide
a description of
$D_4(E,N) \cap N$, also in terms of the exterior
 torsion square; from
this description we deduce an alternate proof of the well-known result that the fourth dimension quotient is always of exponent $\leq 2$. 
\para
Given a group $G$ and an endo-functor $F$ on the category of groups, one can define the derived functors $\mathcal L_iF(G)$ of $F$ to be the homotopy groups $\pi_i(F(\mathfrak F)),$ where $\mathfrak F\to G$ is a free simplicial resolution of $G$. When applied to the functors $\gamma_i/\gamma_{i+1},$  $\gamma_1/\gamma_i,$ $\mathfrak g^i/\mathfrak g^{i+1}$ and $\mathfrak g/\mathfrak g^{i},$ there arise two spectral sequences $E(G)$ and $\overline E(G)$ and a natural map $\kappa: E(G)\to \overline E(G)$. The study of dimension quotients is naturally linked to these spectral sequences; this is explained in section \ref{ss}. For Abelian categories, a theory of non-additive derived functors has been given by Dold-Puppe \cite{DP:61}. Given an Abelian group $A$ and an endo-functor $F$ on the category of Abelian groups, let  $\mathcal D_iF(A)$ denote the derived functor $L_iF(A,\,0)$ in the sense of Dold-Puppe \cite{DP:61}. For every $n\geq 3$, we construct a canonical homomorphism $v_n$ (see diagram (\ref{nachd1})) from a subgroup $V_{n-1}(G)$ of the torsion group $\mathcal D_1S_{n-1}(G_{ab})$ to a sub-quotient of the $n$-th dimension quotient, where $S_i(A)$ denotes the cokernel of the natural map from the $i$-th component $L_i(A)$ of the free Lie ring $L(A)$ on the Abelian group $A$ over $\mathbb Z$ to the $i$-th tensor power of $A$. For this purpose, we need a description of the  derived functors of certain quadratic functors (recalled in section \ref{polfunct})  studied by A. K. Bousfield \cite{Bousfield:67} and  H.-J. Baues and T. Pirashvili \cite{BP:00}. 
After proving an identification theorem in section \ref{it}, we analyse in sections \ref{sew} and \ref{fifthdq} the maps $v_4$ and $v_5$ and provide sufficient conditions for the triviality of the fourth and fifth dimension quotients (see Corollary \ref{4dim} and Theorem \ref{5dim}).  
  \para
For basic notions and properties of simplicial objects, we refer the reader to \cite{May:67}.

\Para
\section{Polynomial functors}\label{polfunct}
\par\vspace{.5cm}
Recall that a map $f:G\to A$ from a multiplicative group $G$ to an Abelian group $A$ is said to be a {\it polynomial map of degree} $\leq n$ if its linear extension to the integral group ring $\mathbb Z[G]$ of $G$ vanishes on $\mathfrak g^{n+1}$, where $\mathfrak g$ denotes the augmentation ideal of $\mathbb Z[G]$ \cite{Passi:68a}.
\par\vspace{.5cm}
Let ${\sf Ab}$ denote the category of Abelian groups. Let $F:{\sf
  Ab}\to {\sf Ab}$ be a functor with $F(0)=0$, where $0$ is the zero
homomorphism. If $M,\ N$ are Abelian groups, then, by definition,
$F$ defines a map $$F_{MN}:Hom(M,\,N)\to Hom(F(M),\,F(N)).$$ Let
us write $U:=Hom(M,\,N)$ multiplicatively, and
$V:=Hom(F(M),\,F(N))$ additively. Then $F_{MN}:U\to V$ is a map
satisfying $F_{MN}(e)=0$, where $e$ is the identity element of
$U$. The functor $F$ is said to be  a {\it polynomial functor of
degree} $\leq n$ if $F_{MN}$ is a polynomial map of degree $\leq
n$ for all Abelian groups $M,\ N$. This is a reformulation
\cite{Passi:69} of the concept of functors of finite degree in the
sense of Eilenberg-MacLane \cite{EM:54a}. 
\Para {\bf Quadratic
functors.}  A functor $F :{\sf Ab}\to {\sf Ab}$ is quadratic, i.e., has
degree $\leq 2$, if $F(0) = 0$  and if the cross effect
$$F(A | B) = \text{Ker}(F(A \oplus B)\to  F(A) \oplus F(B))$$
is biadditive. This yields the binatural isomorphism
$$F(A \oplus B) = F(A) \oplus F(B) \oplus F(A | B)$$
given by $(F_{i_1}; \,F_{i_2};\, i_{12})$ where $i_1 :A\subset A\oplus  B;\  i_2 : A\subset A\oplus  B$ and $i_{12} : F(A|B)\subset F(A\oplus  B)$
are the inclusions. Moreover, for any $ A\in {\sf Ab}$, one gets the diagram \cite{BP:00}\begin{equation}\label{folding}
F\{A\} := (F(A)\buildrel{H}\over \rightarrow
  F(A | A)\buildrel{P}\over \rightarrow
  F(A)).\end{equation}
Here $P = F(p_1 + p_2)i_{12} : F(A|A)\subset F(A \oplus A)\to  F(A)$ is given by the codiagonal $p_1+
p_2 : A\oplus A\to  A$, where $p_1$ and $p_2$ are the projections. Moreover, $H$ is determined by
the equation $i_{12}H = F(i_1+i_2)-F(i_1)-F(i_2)$, where $i_1+i_2 : A \to  A\oplus A$ is the diagonal map.

\para
We  recall from \cite{EM:54a} the definitions of certain  quadratic functors.
\par\vspace{.5cm}
\noindent $\mathbf {\text{\bf Tor}(A,\,C).}$  For Abelian groups $A$ and $C$, the Abelian group
$\text{Tor}(A,\,C)$ (\cite{EM:54a}, \S 11, p.\,85) has generators
$$
(a,\,m,\,c),\ a\in A,\ c\in C,\ 0<m\in \mathbb Z,\ ma=mc=0
$$
and relations
\begin{align*}
& (a_1+a_2,\,m,\,c)=(a_1,\,m,\,c)+(a_2,\,m,\,c),\ \text{if}\ ma_1=ma_2=mc=0\\
& (a,\,m,\,c_1+c_2)=(a,\,m,\,c_1)+(a,\,m,\,c_2),\ \text{if}\ ma=mc_1=mc_2=0\\
& (a,\,mn,\,c)=(na,\,m,\,c),\ \text{if}\ mna=mc=0\\
& (a,\,mn,\,c)=(a,\,m,\,nc),\ \text{if}\ ma=mnc=0.
\end{align*}

In particular, we have the functor $A\mapsto \text{Tor}(A,\,A)$ on
the category {\sf Ab}.
We denote the class of the triple $(a,\,m,\,c)$ by $\tau_m(a,\,c)$.
\para\noindent
$\mathbf{\Omega(A).}$  Let $A$ be an Abelian group. Then the group $\Omega(A)$, defined by Eilenberg-MacLane (\cite{EM:54a}, p.\,93), is the Abelian group
generated by symbols $w_n(x),\ 0<n\in \mathbb Z, x\in A,\ nx=0$ with defining relations
\begin{align*}
& w_{nk}(x)=kw_n(x),\ nx=0,\\
& kw_{nk}(x)=w_n(kx),\ nkx=0,\\
& w_n(kx+y)-w_n(kx)-w_n(y)=w_{nk}(x+y)-w_{nk}(x)-w_{nk}(y),\
nkx=ny=0,\\
& w_n(x+y+z)-w_n(x+y)-w_n(x+z)-w_n(y+z)+w_n(x)+w_n(y)+w_n(z)=0,\\
& nx=ny=nz=0.
\end{align*}
We continue to denote the class of the element $w_n(x)$ in $\Omega(A)$  by $w_n(x)$ itself.
\para
Eilenberg-MacLane (\cite{EM:54a}, p.\,94) constructed a map
$$
E: \text{Tor}(A,\,A)\to \Omega(A)$$by setting
$$
\tau_n(a,\,c)\mapsto w_n(a+c)-w_n(a)-w_n(c).
$$A natural map
\begin{equation}\label{eilmac1}
T: \Omega(A)\to \text{Tor}(A,\,A)
\end{equation}
can be defined by setting
$$
w_n(x)\mapsto \tau_n(x,\,x),\ x\in A,\ nx=0.
$$Clearly the composite map
$$
E \circ T: \Omega(A)\to \Omega(A)
$$
is  multiplication by 2; for, as a consequence of the defining
relations the elements $w_n(x)$ satisfy $$w_n(mx)=m^2w_n(x)$$ for
all $m\in \mathbb Z,\ 0<n\in \mathbb Z,\ x\in A$.
\para\noindent
{\bf Whitehead functor} $\mathbf{\Gamma(A).}$\quad
For $A\in {\sf Ab}$, define the commutative graded ring $\mathbf{\Gamma}(A)$
generated by the elements $\gamma_t(x)$ for each $x\in A$ and each
non-negative integer $t$, of degree $2t$, subject to the relations
\begin{align}
\gamma_0(x)&=1,\\
\gamma_t(rx)&=r^t\gamma_t(x),\\
\gamma_t(x+y)&=\sum_{i+j=t}\gamma_i(x)\gamma_j(y),\\
\gamma_s(x)\gamma_t(x)&=\binom{s+t}{s}\gamma_{s+t}(x).
\end{align}

The homogeneous component $\mathbf{\Gamma}_4(A)$ of $\mathbf{\Gamma}(A)$ of degree 4 can be identified with the Whitehead  functor (\cite{Whitehead:50},
\cite{EM:54a}, pp.\,92 \& \,110)
$\Gamma:{\sf Ab}\to {\sf Ab},\ A\mapsto \Gamma(A)$, where the group $\Gamma(A)$ is defined for $A\in {\sf Ab}$ to be
the group  given by generators $\gamma(a)$, one for each $x\in A$,
subject to the defining relations\begin{equation}
\gamma(-x)=\gamma(x),
\end{equation}
\begin{equation}
\gamma(x+y+z)-\gamma(x+y)-\gamma(x+z)-\gamma(y+z)+\gamma(x)+\gamma(y)+\gamma(z)=0
\end{equation}
for all $x,\,y,\,z\in A$.
\para\noindent
$\mathbf{R(A).}$  For $A\in {\sf Ab}$, let $_2A$ denote the subgroup consisting of
elements $x$ satisfying $2x=0$. Define  $R(A)$ (\cite{EM:54a}. p.\,120)
to be the quotient group of
Tor$(A,\,A)\oplus \varGamma(_2A)$ by the relations
\begin{align}
\tau_m(x,\,x)&=0,\quad mx=0,\\
\gamma_2(s+t)-\gamma_2(s)-\gamma_2(t)&=\tau_2(s,\,t),\quad s,\,t\in _2A.
\end{align}
\para
The functors $A\mapsto \Gamma(A),\ \text{Tor}(A,\,A),\ \Omega(A),\
R(A)$ are all quadratic functors on the category of Abelian
groups, i.e., all these functors have  degree $\leq 2$.
Furthermore, $R(A) = H_5K(A; 2);\ \Omega(A)= H_7K(A;
3)/(Z/3Z\otimes A)$, $ R(A | B) = \Omega(A| B) =\text{Tor}(A,\,B)$
and $R(\mathbb Z) = \Omega(\mathbb Z) = 0$ and $R(\mathbb Z/n) =
\mathbb Z/(2,\,  n);\ \Omega\mathbb (Z/n) = \mathbb Z/n$. \Para
{\bf The  square functor.} Let  ${\sf Ab}_g$ denote  the category
of graded Abelian groups. For  $A,\ B\in {\sf Ab}$,  let $A\otimes
B$ and $A*B:=\text{Tor}(A,\,B)$ be respectively the tensor product
and the torsion product of $A,\ B$.  The notion of tensor product
of Abelian groups extends naturally to that of {\it tensor product
$A\otimes B$ of graded Abelian groups} $A,\ B$ by setting
$$
(A\otimes B)_n=\bigoplus_{i+j=n}A_i\otimes B_j.
$$
We also need the {\it ordered tensor product}
$A\buildrel{>}\over\otimes B$ of graded Abelian groups, which is defined by
setting $$
(A\buildrel{>}\over{\otimes} B)_n=\bigoplus_{i+j=n,\
i>j}A_i\otimes B_j
$$for $A,\ B\in {\sf Ab}_g$.
In an analogous manner, we can define, for $A,\ B\in {\sf Ab}_g$, {\it torsion product} $A*B$ and  {\it ordered torsion product}
$A\buildrel{>}\over{*}B$ as
$$
A*B=\bigoplus_{i+j=n}A_i* B_j,\quad
(A\buildrel{>}\over{*}B)_n=\bigoplus_{i+j=n,\ i>j}A_i* B_j.
$$

The tensor product, torsion product and the ordered tensor and
torsion product are, in an  obvious way, bifunctors on the category {\sf Ab}$_g$.

\par\vspace{.5cm} Let $\wedge^2$ be  the exterior square functor on
the category {\sf Ab}. The {\it weak square functor}
$$
sq^\otimes: {\sf Ab}_g\to {\sf Ab}_g
$$ is defined  by
$$
sq^\otimes(A)_n=\begin{cases} \Gamma(A_m),\ \text{if}\ n=2m,\ m\ \text{odd},\\
\wedge^2(A_m),\ \text{if}\ n=2m,\ m\ \text{even},\\
0,\ \text{otherwise}.\end{cases}
$$
Let $(\mathbb Z_2)_{odd}$ be the graded Abelian group which is
$\mathbb Z_2$ in odd degree $\geq 1$ and  trivial
otherwise; thus $(\mathbb Z_2)_{odd}$ is the reduced homology of
the classifying space $\mathbb RP_\infty=K(\mathbb Z_2,1)$. The  {\it square functor} $Sq^\otimes: {\sf Ab}_g\to {\sf Ab}_g$ is  defined as follows:
$$
Sq^\otimes(A)=A\buildrel{>}\over{\otimes}(A\oplus (\mathbb
Z_2)_{odd})\oplus sq^\otimes(A).
$$Clearly the square functor is a quadratic functor, its second
cross-effect functor (\cite{EM:54a}, p.\,77) is
$$
Sq^\otimes(A|B)=A\otimes B,$$and one has the operators
\begin{equation}
Sq^\otimes(A)\buildrel{H}\over{\rightarrow} A\otimes
A\buildrel{P}\over{\rightarrow} Sq^\otimes(A)
\end{equation}
which are induced by the diagonal and the folding map
respectively [see (\ref{folding})].
\par\vspace{.25cm}
Define next  the {\it torsion square functor} $$Sq^\star(A):{\sf
Ab}_g\to {\sf Ab}_g$$ by setting
$$
Sq^\star(A)=(A\buildrel{>}\over{*}(A\oplus (\mathbb
Z_2)_{odd}))\oplus sq^\star(A),
$$
where
$$
sq^\star(A)_n=
\begin{cases} \Omega(A_m),\ n=2m,\ m\ \text{even}\\
R(A_m),\ n=2m,\ m\ \text{odd}\\
0,\ \text{otherwise}
\end{cases}
$$
\par\vspace{.5cm}
\section{Relative dimension subgroups and polynomial (co)homology}\label{relative1}
\para
Relative dimension subgroups naturally appeared in the study (\cite{Passi:68b}, \cite{Passi:79}) of the dimension subgroup problem
via polynomial cohomology  introduced there for this purpose. For any normal subgroup $N$ of a group $E$, define the $n$-th dimension subgroup $D_n(E,\,N)$ of $E$ relative to $N$ by setting \[ D_n(E,\,N) = E \cap \Big( 1 + \mathfrak n\mathfrak e + \mathfrak e^{n} \Big)\;.\]
This group is also called {\em dimension subgroup of the extension}\/ $1 \to N \to E \to E/N\to 1$ by Kuz'min who exhibits many interesting properties in \cite{Kuzmin:96}. Note that $D_n(E) = D_n(E,\,\{1\}) = D_n(E,\,\gamma_{n-1}(E))$ since $\gamma_{n-1}(E)-1 \subset \mathfrak e^{n-1}$, whence relative dimension subgroups generalize  ordinary dimension subgroups; and indeed, they were extensively used to study the latter since the work in  \cite{Passi:68b}. 
\para
Note that $N'\gamma_n(E) \subset D_n(E,N)$, where $N'$ denotes  the derived subgroup $\gamma_2(N)$ of $N$, so the {\em relative dimension subgroup problem}\/ asks under which conditions equality holds.
As $D_2(E,\,N) = D_2(E) = \gamma_2(E)$,   the problem becomes interesting only for $n\ge 3$. The basic tool for its study consists of the investigation of the following  constructions. Let
\[ P_n(E,\,N) = \mathfrak e/(\mathfrak n \mathfrak e + \mathfrak e^{n+1}),\quad I(E,\,N)=\mathfrak e/\mathfrak n\mathfrak e, \] 
\[ P_n(E) = \mathfrak e/\mathfrak e^{n+1} = P_n(E,\,\{1\}) = P_n(E,\,\gamma_{n}(E)).\] 
Let $G=E/N$. Then the group $P_n(E,\,N)$ has the structure of both a non-unitary ring and a left nilpotent $G$-module of class $\le n$, via left multiplication in $E$, and the map $$p_n\mapco E\to P_n(E,\,N),\quad p_n(e) =  \overline{e-1},$$ is a $\pi $-derivation, where  $\pi\mapco E\auf G$ is the natural projection, i.e. $p_n$ satisfies $$p_n(ee')= \pi(e) p_n(e') + p_n(e) = 
 p_n(e ) + p_n(e') +  p_n(e)  p_n(e')\  \text{for}\ e,\ e'\in E.$$ Recall the   classical sequence of $G$-modules
 \BE\label{HiStsequ}  0 \to N/N' \mr{d}I(E,\,N) \mr{I(\pi)} \mathfrak g \to 0\EE
 with  $d(\bar{n})=\overline{n-1}$  and $I(\pi)(\,\overline{e-1}\,) = \pi(e)-1$ (see \cite{HS:71}); reducing modulo $\mathfrak e^{n+1}$ induces a sequence of nilpotent $G$-modules of class $\le n$
 \BE\label{Pnsequ} 1 \to \frac{\dst D_{n}(E,\,N) \cap  N\gamma_n(E)}{\dst N'\gamma_n(E)} \hmr{}
  \frac{\dst    N\gamma_n(E)}{\dst N'\gamma_n(E)}
 \hmr{\bar{d}} P_{n-1}(E,\,N) \mr{P_{n-1}(\pi)}
 P_{n-1}(G) \to 0
 \EE
This sequence then implies a short exact sequence
\BE\label{Dnrelsequ}
1  \hmr{} \frac{\dst D_{n}(E,\,N) \cap  N\gamma_n(E)}{\dst N'\gamma_n(E)} \hmr{}
\frac{\dst D_{n}(E,\,N) }{\dst N'\gamma_n(E)}
\hmr{} \frac{\dst D_n(G)}{\dst \gamma_n(G)}  \hmr{} 1.
\EE
Consequently,  the study of relative dimension subgroups breaks down into the study of the ``intersection quotient"
$$ \frac{ D_{n}(E,\,N) \cap  N\gamma_n(E)}{ N'\gamma_n(E)}$$ and of ordinary dimension subgroups of groups of lower nilpotency class; indeed, $N$ can be replaced by $N\gamma_{n-1}(E)$ since $D_{n}(E,\,N)= D_{n}(E,\,N\gamma_{n-1}(E))$.
Thus one may focus on the study of the intersection quotient.
\para
It turns out that there is a qualitative difference between the case where $N$ is central (or $[E,N] \subset N'\gamma_n(E)$) and the non-central case. Indeed, for $n=3$ the intersection quotient (which actually equals $D_3(E,N)/N'\gamma_3(E)$) was proved to be trivial if $N$ is trivial or $G$ is cyclic \cite{Passi:79}, and under some other, notably torsion-freeness conditions \cite{KV:88}; but it was completely determined - and shown to be non-trivial in general - only by Hartl in an unpublished note in 1993. We cite the result here for further use; its proof, along with a generalization to arbitrary coefficient rings and to relative  Fox subgroups, can be found in \cite{Hartl:D3F2}.
\para
Consider the following part of a six-term exact sequence for the tensor and torsion product of
Abelian groups (see \cite[V.6]{Maclane:63}).
\[ {\rm Tor}(G\ab,\,G\ab) \hmr{\omega_1} (NE'/E') \ot G\ab 
\hmr{i \ot1} E\ab \ot G\ab \hmr{q \ot 1} G\ab \ot G\ab \hmr{} 0 \]
 Here $i,\ q$ are the canonical injection and
projection, respectively. Moreover, commutation in $E$ induces a homomorphism
 \[ \map{[\,,\,]}{ (NE'/E') \ot G\ab }{N\gamma_3(E)/N'\gamma_3(E)}\,.\]
\par\vspace{.25cm}
\begin{satz}\label{D3rel}\quad For any normal subgroup $N$  of a group $E$  the following
sequence of natural homomorphisms is exact:
  \[ {\rm Tor}_1^{\Z}(G\ab,\, G\ab) \hmr{[\,,\,] \omega_1} N\gamma_3(E)/N'\gamma_3(E) \hmr{\bar{d}}
P_2(E,N) \hmr{P_2(\pi)} P_2(G) \hmr{} 0. \]
Consequently,  one has
\[ D_3(E,\,N) = N'\gamma_3(E)\, sgr\{[a^k,\,b]\,|\,\mbox{$a,\,b\in E$, $k\in \Z$, $a^k,b^k\in NE'$}\} \,.\]
\end{satz}
\par\vspace{.25cm}
We note that the formula for $D_3(E,N)$ provided by the above theorem was also reproved by different methods in \cite{RTV:98}.\V

Theorem \ref{D3rel} surprisingly shows that unlike the case where $N$ is central, the inclusion 
$N'\gamma_3(E) \subset D_3(E,\,N) $ is   {\it not  always} an equality, as was
suggested by the known partial results. Indeed, there are counter-examples which
are $p$-groups for {\em any prime} $p$, see \ref{cex} below; this is in contrast to the case of
classical dimension subgroups (i.e.\ $N=1$) which coincide with the terms of the lower central
series of $G$ unless $p=2$ (cf.\ \cite{Gupta:02}).
\para
\begin{bsp}\label{cex}\quad \rm Let $p$ be a prime and $0<r\le s$.
Define
  \[ E=\langle x,\,y\, |\,1= x^{p^{s+1}} = y^{p^{s+1}} = [x,\,[x,\,y]] = [y,\,[x,\,y]]\,\rangle \:.\]
Let $N = sgp\{x^{p^r},\, y^{p^{s}},\, [x,\,y]\}$. Then $z = [x,\,y]^{p^{s}} =
[x,\,y^{p^{s}}] \in D_3(E,\,N)$ by Theorem \ref{D3rel}, but $z$ has order $p$ modulo
$N'\gamma_3(E) = \{1\}$.\hfbox

\end{bsp}
\para
Observe that  Theorem \ref{D3rel}   implies that the group $D_3(E,N)/N'\gamma_3(E)$ is a quotient of $D_3(F,\,R)/R'\gamma_3(F)$ for any free presentation $1\to R\to F\to G \to 1$ of $G$, and that the latter quotient  does not depend on the chosen presentation. We will now show that it even is a functor in $G\ab$ which identifies with another functor which can be easily described as follows.

\para
Let {\sf Ab} denote the category of Abelian groups, and the
 functor $\stackrel{\sst\sm}{*}\mapco {\sf Ab} \to {\sf Ab}$ be defined by
\BE\label{LSdef} A\stackrel{\sst\sm}{*}A = {\rm Tor}(A,\,A) / sgr\{\tau_m(x,\,x)\,|\,x\in A\}  \EE
which may be called the {\em exterior torsion square}\/ of $A$; this functor is, in fact,  the first derived functor (in the sense of \cite{DP:61}) of the symmetric square (see section \ref{ss}).
\para
\begin{satz}\label{D3free} Let $R$ be a normal subgroup of a free group $F$. Then the map $[\,,\,] \omega_1$ in Theorem \ref{D3rel} induces a natural isomorphism
\[   (F/RF')\stackrel{\sst\sm}{*} (F/RF')\hcong D_3(F,\,R)/R'\gamma_3(F) \,.\]
\end{satz}
\para
The proof requires the functor $\Gamma\mapco {\sf Ab} \to {\sf Ab}$ (see section \ref{polfunct}) and the natural homomorphisms
\[ A\ot A \mr{w} \Gamma(A) \mr{\delta} A\ot A \Sur{q_{\sm}} \Lambda^2(A) = \coker{\delta} \]
where $w(x\ot y) = \gamma(x+y)-\gamma(x)-\gamma(y) = w(y\ot x)$ and $\delta\gamma(x)=x\ot x$ for $x,\,y\in A$. Note that $\delta w(x\ot y) = x\ot y +y\ot x$. We also need the standard fact that, for an exact sequence $A\mr{f}B\mr{g} C\to 0$
of Abelian groups, the induced sequence \[\Gamma(A) \oplus A\ot B \mr{\alpha} 
\Gamma(B) \mr{\Gamma(g)}  \Gamma(C) \to 0 \]
is exact for $\alpha = (\Gamma(f)\,,w(f\ot 1))$.\vspace{5mm}

\def\VC#1#2#3{\makebox[0cm][l]{\hspace*{-#1mm}\makebox[#2cm]{\raisebox{4mm}{$\scriptstyle#3$}}}
\begin{minipage}{0cm} \unitlength1cm \hspace*{-#1mm}
\begin{picture}(1,0) \put(0,0){\vector(1,0){#2}}\end{picture}  \end{minipage}
}

\N{\bf Proof of Theorem \ref{D3free}\,:\hspace{2mm}} According to Theorem \ref{D3rel}, we must show that $$\Ker{[,]\omega_1}  = sgr\{\tau_m(x,\,x)\,|\,x\in A\} \ \text{for}\ A=F/RF'.$$
Abbreviating  $D=F\ab$ and $C=RF'/F'$ we get a \Z-free resolution  $$C\stackrel{i}{\hra} D \Sur{q} A.$$ Note that we have an injection and an  isomorphism $$\frac{R\cap\gamma_2(F)}{R'\gamma_3(F)} \stackrel{j}{\hra}
\frac{ \gamma_2(F)}{R'\gamma_3(F)} \mr{\bar{c}^{-1}} \frac{\Lambda^2(D)} {\Lambda^2(C)},$$ where $$c\mapco \Lambda^2(D) \isor{}  \gamma_2(F)/\gamma_3(F)$$ is the classical  isomorphism of Witt in degree 2, i.e., $c((aF')\sm (bF')) = [a,b]\gamma_3(F)$ for $a,b\in F$. Thus $\Ker{[,]\omega_1} = 
\Ker{\beta}$ for $\beta =\bar{c}^{-1} j [,]\omega_1$.
Now consider the following commutative diagram
\[\begin{matrix}
{\rm Tor}(A,\,A)  & \mr{\beta}  &  \frac{\dst \Lambda^2(D)} {\dst \ruled \Lambda^2(C)} & & & & & \cr
\mapdown{\omega_1} & & \surup{\overline{q_{\Lambda}}} & & & & & \cr
\ruled C\ot A & \mr{\overline{i \ot 1}} & \frac{\ruleu\dst D\ot D}{\ruled \dst C\ot C} & \mr{\overline{q \ot 1}} & A\ot D & \mr{1\ot q} & A\ot A \cr
 & & \mapup{\bar{\delta}} & & & & \mapup{\delta}\cr
\Gamma(C) \oplus A\ot C & \mr{\alpha} & \frac{\ruleu\dst \Gamma(D)}{\dst \ruled w(C\ot C)}  &\VC{2}{3}{\overline{\Gamma(q)}} & & &  \Gamma(A)
  \end{matrix}\]
whose last row is exact with $$\alpha= (q_{\Gamma}\Gamma(i),\, \overline{w(1\ot i)})$$ and $q_{\Gamma}\mapco \Gamma(D) \auf \Gamma(D)/w(C\ot C)$ being the canonical projection. As $\Imm{\,\overline{i \ot 1}\,} = \Ker{\,\overline{q \ot 1}\,}$,  we get
\begin{eqnarray}
\Ker{[,]\omega_1}  &=& \omega_1^{-1} \, (\overline{i \ot 1} \, )^{-1}  \, \Ker{\overline{q_{\Lambda}}} \nonumber\\
 &=& \omega_1^{-1} \, (\overline{i \ot 1} \, )^{-1}  \, \Imm{\bar{\delta}} \nonumber\\
 &=&\omega_1^{-1} \, (\overline{i \ot 1} \, )^{-1}  \,  \bar{\delta}  \, \Ker{(\overline{q \ot 1} \,)\bar{\delta} }\label{Kerred}
 \end{eqnarray}
 Note that $\omega_1$ and $\overline{i \ot 1}$ are injective as $D$ is a free \Z-module, so we essentially must determine the group $\Ker{(\overline{q \ot 1} \,)\bar{\delta} }$. To this end, consider the following commutative diagram with exact rows and columns, where $\omega_2$ is the corresponding connecting homomorphism; moreover, $\alpha'$ and $\overline{\Gamma(q)}'$ are given by restriction of $\alpha$ and $\overline{\Gamma(q)}$:
\[\begin{matrix}
\ruled & & \Gamma(C) & \mr{\alpha'} & \Ker{(\overline{q \ot 1} \,)\bar{\delta} } & \mr{\overline{\Gamma(q)}'}  & \Ker{\delta} & & \cr
& &\mapdown{(1,\,0)^t} & & \injdown{} & &  \injdown{}\cr
\ruled & & \Gamma(C)  \oplus A\ot C & \mr{\alpha} &
\frac{\dst\ruleu  \Gamma(D)}{\dst \ruled w(C\ot C)}  &\mr{\overline{\Gamma(q)}} &  \Gamma(A) & \to & 0\cr
& & \surdown{(0,\,1)}  &  & \mapdown{(\overline{q \ot 1} \,)\bar{\delta} } & & \mapdown{\delta} \cr
\ruleu{\rm Tor}(A,\,A) & \mr{\omega_2} & A \ot C & \mr{1\ot i} & A\ot D & \mr{1\ot q} & A\ot A & \to & 0
  \end{matrix}\]\par\noindent
  Applying the snake lemma provides an exact sequence
  \BE\label{Kersequ} \Gamma(C)  \oplus {\rm Tor}(A,\,A) \hmr{\alpha''} \Ker{(\overline{q \ot 1} \,)\bar{\delta} } \hmr{\overline{\Gamma(q)}'}    \Ker{\delta} \to 0, \EE
  where $\alpha'' = (q_{\Gamma}\Gamma(i)\,, \overline{w(1\ot i)}\omega_2)$. Now for $c\in C$,
  \[\bar{\delta} \alpha'' \gamma(c) = 
  \bar{\delta} q_{\Gamma}\Gamma(i)\gamma(c) = 
  \bar{\delta}  \, \overline{ \gamma(c)} =  
  \overline{\delta \gamma(c)} = \overline{c\ot c} = 0\,.\]
  Next, let $ \tau_m(x_1,\,x_2) \in {\rm Tor}(A,\,A)$, and let $x_k'\in D$ such that $qx_k' = x_k$, $k=1,\,2$. Then 
  \begin{eqnarray*}
\bar{\delta} \alpha'' \tau_m(x_1,\,x_2) &=& \bar{\delta}\,\overline{w(1\ot i)\,}  \omega_2 \tau_m(x_1,\,x_2) \\
&=& \bar{\delta} \,\overline{w(1\ot i)} (x_1\ot i^{-1}(mx_2')) \\
 &=& \delta w (x_1' \ot mx_2') + C\ot C \\
 &=& x_1' \ot mx_2' + mx_2' \ot  x_1' + C\ot C \\
 &=& (i\ot 1) ( i^{-1}(mx_1') \ot x_2' + i^{-1}(mx_2') \ot x_1' ) + C\ot C \\
  &=& (\,\overline{i\ot 1}\,) \omega_1 (  \tau_m(x_1,\,x_2) + \tau_m(x_2,\,x_1) )
\end{eqnarray*}
Thus   $\omega_1^{-1} \,(\,\overline{i\ot 1}\,)^{-1} \bar{\delta} \,\Imm{\alpha''}$ is generated by the elements 
  $\tau_m(x_1,\,x_2) + \tau_m(x_2,\,x_1)$, $x_1,\,x_2\in A$ which belong to  $sgr\{\tau_m(x,\,x)\,|\,x\in A\} $.
  \para
  By  \REF{Kersequ} it remains to treat \Ker{\delta}. We know that it is generated by the elements $o(x)\gamma(x)$, where $x\in A$ is of finite even order $o(x)$ (see \cite[Lemma 5.1]{Hartl:96}). If $x$ is of odd order, then 
  $o(x)\gamma(x)=0$; in fact, $w(x\ot x) = \gamma(2x) -2\gamma(x) = 2\gamma(x)$, whence $o(x)\gamma(x)=
  o(x)\gamma(x) + \binom{o(x)}{2} w(x\ot x) = o(x)^2\gamma(x) = \gamma(o(x)x) = 0$. Thus \Ker{\delta} is generated by the elements $m\gamma(x)$ where $x\in A$, $m\in \Z$ such that $mx=0$. For such $x$ and $m$, and $x'\in D$ such that $qx'=x$, we have
    \begin{eqnarray*}
\bar{\delta} q_{\Gamma} (m\gamma(x') ) &=& m \,\overline{\delta \gamma(x')} \\
&=& (i\ot 1) ( i^{-1}(mx') \ot x') + C\ot C\\
 &=&  (\,\overline{i\ot 1}\,) \omega_1    \tau_m(x ,\,x ) \,.
  \end{eqnarray*}
Thus it follows from equation \REF{Kerred} that 
$\Ker{[,]\omega_1} = sgr\{\tau_m(x,\,x)\,|\,x\in A\} $, as asserted.\hfbox\V

Theorem \ref{D3free}  admits an application to group cohomology $H^2_{{\mathcal N}_c}(G,\,M)$ with respect to the variety ${\mathcal N}_c$
of $c$-step nilpotent groups, for $c=2$, thus solving a problem of Leedham--Green \cite{LG:71} in this case.  Let $G$ be nilpotent of class $c$ and $M$ be a nilpotent $G$-module of class $\le c$, i.e., $\mathfrak g^{c} M=0$. Then it is well-known that the set of congruence classes of group extensions  $e\mapco0\to M \mr{i}E \mr{\pi} G \to 1$ with $E\in {\mathcal N}_c$, equipped with the Baer sum operation, is a group isomorphic with $H^2_{{\mathcal N}_c}(G,\,M)$.
 Leedham--Green asks when the usual shifting isomorphism ${\rm Ext}^1_{\Z[G]}(\mathfrak g,\,M) \cong H^2(G,\,M)$ admits a varietal analogue which for ${\mathcal N}_c$ means a natural isomorphism ${\rm Ext}^1_{R_c}(P_c(G),\,M)
  \cong H^2_{{\mathcal N}_c}(G,\,M)$ with $R_c = \Z[G]/\mathfrak g^{c}$. He shows that even for $c=2$ this is not the case, by taking $G$ to be the Klein four-group. Kuz'min \cite{Kuzmin:83} establishes a natural exact sequence which, for the variety ${\mathcal N}_c$ and in our notation, reads
  \BE\label{Kusequ}
0 \hmr{}  {\rm Ext}^1_{R_c}(P_c(G)\,,M) \mr{\mu} H^2_{{\mathcal N}_c}(G,M) \mr{\nu} {\rm Hom}_{R_c}(
 (D_{c+1}(F,R) \cap R)/ R',M) \EE
where $f\mapco 1 \to R \to F \mr{\varphi} G \to 1$ is a free presentation of $G$ {\em in the variety ${\mathcal N}_c$}\/, and where the structure of
 $R_c$-module on
$(D_{c+1}(F,R) \cap R/R'$  is given by the usual
 conjugation action of
$G$ on $R/R'$.  Moreover, $\nu$ is given as follows: for an extension $e$ in ${\mathcal N}_c$ as above choose a lifting $\alpha_0$ giving rise to a commutative diagram
\[\begin{matrix}
\ruled1 &\to& R/R'  & \to & F/R' & \to & G & \to &1\cr 
 & & \mapdown{\alpha_1} & &\mapdown{\alpha_0} & & \| & &\cr
\ruleu 1 & \to & M & \mr{i} & E &\mr{\pi} & G & \to & 1
\end{matrix}\]
Then $\nu[e]$ is given by restriction of $\alpha_1$. For $c=2$, we improve Kuz'min's sequence as follows.\para
 Consider the group $G\ab \stackrel{\sst\sm}{*} G\ab$ as a trivial $R_c$-module.
By Theorem \ref{D3free} and \REF{Pnsequ}, we have two successive short exact sequences of $R_c$-modules
 \BE\label{2sequs} G\ab \stackrel{\sst\sm}{*} G\ab \hspace{1mm}\Inj{\overline{[,]\omega_1}} \hspace{1mm}R/R' \hspace{1mm}\Sur{}\hspace{1mm} \coker{[,]\omega_1}\hspace{1mm}\Inj{\bar{d}} \hspace{1mm}P_2(F,\,R) \Sur{P_2(\varphi)} P_2(G).\EE
 These sequences give rise to two successive boundary maps\V
 
 \N$\hfill {\rm Hom}_{R_2}(G\ab \stackrel{\sst\sm}{*} G\ab,\,M) \hmr{\partial_1}
{\rm Ext}_{R_2}^1( \coker{[,]\omega_1},\,M) \hmr{\partial_2}
{\rm Ext}_{R_2}^2(P_2(G),\,M). \hfill$\vspace{5mm}


\begin{satz}\label{Kuc=2} If $G$ is a $2$-step nilpotent group and $M$  nilpotent $G$-module of class $\leq 2$, then   there is a natural exact sequence 
\begin{multline*}
 0 \to {\rm Ext}^1_{R_2}(P_2(G),\,M) \mr{\mu} H^2_{{\mathcal N}_2}(G,\,M) \mr{\nu'}  {\rm Hom}_{R_2}(G\ab \stackrel{\sst\sm}{*} G\ab,\,M)\\ \mr{\partial} {\rm Ext}^2_{R_2}(P_2(G),\,M), \end{multline*}
where $ \partial  =  \partial_2\partial_1 $, and $\nu'$ is
defined as follows:
\para
 For a typical generator 
$\tau_m(x_1,\,x_2)$ of ${\rm Tor}(G\ab,\,G\ab)$, let
$\T{x}_k \in E$ be such that $\pi(\T{x}_k)G'=x_k$, and let $y_2\in M$ be  such that $\T{x}_2^m \equiv i(y_2)$ mod $E'$. Then
$$\nu'[e](\,\overline{\tau_m(x_1,\,x_2)}\,) = (\pi(\T{x}_1)-1)y_2.$$ 
\end{satz}
\par\vspace{.25cm}
We note that the extension of Kuz'min's sequence \REF{Kusequ} on the right by a boundary operator $\partial$ as in the Theorem can be easily generalized to any variety of groups.
\para
\begin{bsp} Let $E$ be the group in Example \ref{cex},  $$M=\Z/p^{s-r+1}\langle x^{p^r} \rangle \times 
\Z/p\langle y^{p^s} \rangle \times 
\Z/p^{s+1}\langle [x,\,y] \rangle   \subset E,$$ and  $G=E/M = \Z/p^{r}\langle \bar{x} \rangle \times 
\Z/p^s\langle \bar{y} \rangle$. Then, for $e\mapco M \hra E \auf G$, one has $[e]\notin \Imm{\mu} $, since $\nu'[e](\,\overline{\tau_{p^r}(\bar{x}, \,p^{s-r}\bar{y})}\,) = p^s[x,\,y] \neq 0$.
\end{bsp}

\par\vspace{.25cm}
\N{\bf Proof of Theorem \ref{Kuc=2}\,:}\hspace{0mm} 
We use the isomorphism $$\overline{[,]\omega_1}\mapco G\ab \stackrel{\sst\sm}{*} G\ab \isor{} D_3(F,\,R)/R'$$ of Theorem \REF{D3free} to replace $\nu$ in \REF{Kusequ} by $$\T{\nu}\mapco [e] \mapsto 
(\overline{[,]\omega_1})^*\nu[e] \in {\rm Hom}_{R_2}(G\ab \stackrel{\sst\sm}{*} G\ab,\,M).$$ But 
choosing  $\hat{x}_k \in F$ to be such that $\varphi(\hat{x}_k)G'=x_k$ and taking $\T{x}_k = \alpha_0(\hat{x}_k)$, $k=1,\,2$, we have
$i_*\T{\nu}[e]\big(\overline{\tau_m(x_1,\,x_2)} \big) = i \alpha_1[\hat{x}_2^m,\,\hat{x}_1]^{-1} = [\alpha_0(\hat{x}_1) ,\,\alpha_0(\hat{x}_2)^m] = [\T{x}_1 ,\,i(y_2)] =
i  \big(  (\pi(\T{x}_1) - 1) y_2) \big) = i_*\nu'[e](\,\overline{\tau_m(x_1,\,x_2)}\,)$. Thus $\T{\nu} = \nu'$ and our sequence is exact at $H^2_{{\mathcal N}_2}(G,M)$. Now let $\beta\in {\rm Hom}_{R_2}(G\ab \stackrel{\sst\sm}{*} G\ab,\,M)$. One has $\partial(\beta) =0$ if and only if  
$\partial_1(\beta) =0$, since  $\partial_2$ is an isomorphism; in fact, if $(x_i)_{i\in I}$ is a basis of the ${\mathcal N}_2$-free group $F$ then $P_2(F,\,R)$ is a free $R_2$-module with basis $(p_2(x_i))_{i\in I}$, as follows from the corresponding fact that the augmentation ideal of an absolutely free group $\Gamma$ is a free $\Z(\Gamma)$-module. 
But $\partial_1(\beta) =0$ if and only if  there exists an $R_2$-linear map $\T{\beta}\mapco R/R' \to M$ such that $
\T{\beta}\,\overline{[,]\omega_1} = \beta$. If such 
$\T{\beta}$ exists then there is an induced extension  $\T{\beta}_*f\mapco 0\to M\to E_{ \T{\beta}} \to G\to 1$ for which we can choose $\alpha_0$ such that $\alpha_1=\T{\beta}$, whence $\nu'[\T{\beta}_*f]=
\T{\nu}[\T{\beta}_*f]=\beta$.
\para
 Conversely, if there exists an extension  
$e\mapco 0\to M\to E  \to G\to 1$ such that $\nu'[e]=\beta$ then the  map $\T{\beta} =\alpha_1$ satisfies $\T{\beta}\,\overline{[,]\omega_1} = \T{\nu}[e] = \nu'[e]= \beta$. Thus our sequence is also exact at ${\rm Hom}_{R_2}(G\ab \stackrel{\sst\sm}{*} G\ab,\,M)$.\hfbox

\para
For $n\ge 4$,  the intersection quotient $ \frac{ D_{n}(E,N) \cap  N\gamma_n(E)}{ N'\gamma_n(E)}$ 
appears to be rather mysterious 
if $N$ is non-central; in fact, nothing seems to be known in this
case. For central $N$, however,   it is known that  the triviality of the intersection quotient  is measured by the polynomial cohomology group $P_nH^2(E/N,\,\Q/\Z)$, where $P_nH^2(G,\,A)$, for a group $G$ and a trivial $G$-module $A$, is defined to be the subgroup of $ H^2(G,A)$ consisting of elements representable by bipolynomial cocyles of degree $\le n$ in each variable. 
\para
In \cite{Hartl:96a} and \cite{Habil}, the above approach   was modified and extended so as to define polynomial cohomology groups $P_nH^i(G,\,M)$ for all $i\ge 1$ and $n$-step nilpotent $G$-modules $M$, and to interpret them as relative Ext-functors. It turned out that polynomial 2-cohomology in this sense provides a powerful tool for {\em computing}\/ the actual cohomology group $H^2(G,M)$ if $G$ is finitely generated torsion-free nilpotent and $M$ is additively torsion-free, not only as a group but along with {\em explicit}\/ representing 2-cocycles in terms of numerical polynomial functions \cite{Hartl:96a}. Concerning relative dimension subgroups a dual version of the  original cohomological approach is proposed in \cite{Diss} and \cite{Habil}, using polynomial homology $P_nH_2(G)$ (instead of cohomology); the reason being that the kernel of the canonical approximation $\rho_{n-2}\mapco H_2(G) \auf P_{n-2}H_2(G)$ turns out to be the {\em universal}\/ relative  dimension quotient $D_n(E,C)\cap N/\gamma_n(E)\cap N$ for all central extensions 
$1 \to C \to E \to G\to 1$, provided  $\gamma_n(G)=1$.\vspace{5mm}

To describe the above  concept, we first recall from \cite{Habil} the construction and basic properties of polynomial (co)homology in the generalized sense.
 \vspace{5mm}

Let $K$ be a commutative unitary ring and $A$ be an augmented $K$-algebra with augmentation ideal $\bar{A}$. 
\para
\begin{defi} \rm  Let $n\ge 0$. The polynomial bar construction  $(P_nB(A)\,,\,\bar{\delta})$, of degree $n$ over $A$, is defined by $ P_nB_0(A) =0$ and 
\[P_nB_i(A) = (A/\bar{A}^{n+1}) \ot_K (\bar{A}/\bar{A}^{n+1})\htt{i-1} \ot_K (\bar{A}/\bar{A}^{n+2}) \]
for $i\ge 1$, and the differential $\bar{\delta}_i\mapco P_nB_i(A) \to P_nB_{i-1}(A)$ is given by
\[ \bar{\delta}_i(\overline{a_0} \ot \cdots \ot \overline{a_i}) = \sum_{j=0}^{i-1} (-1)^j \overline{a_0} \ot \cdots \ot \overline{a_j}   \overline{a_{j+1}} \ot \cdots \ot \overline{a_i} \]
for $i\ge 2$. For   left (resp.\ right) $(A/\bar{A}^{n+1})$-module $M$ (resp. $N$), define {\em polynomial (co)homo\-logy of degree $n$ of A}\/ by
\[ P_nH^i(A,\,M) = H^i({\rm Hom}_{A/\bar{A}^{n+1}}(P_nB(A)\,,\,M)), \]
\[ P_nH_i(A,\,N) = H_i(N \ot_{A/\bar{A}^{n+1}} P_nB(A))\ruled \]
Now, for a group $G$, polynomial (co)homo\-logy of degree $n$ is defined by applying these constructions to the group ring $A=\Z[G]$: for left (resp.\ right) $(n+1)$-step nilpotent $G$-module $M$ (resp. $N$) let
\[P_nH^i(G,\,M) = P_nH^i(\Z[G],\,M), \] \[P_nH_i(G,\,N) = P_nH_i(\Z[G],\,N)\,.\] 
\end{defi}\vspace{3mm}

In order to understand  the properties of these constructions, recall the 
bar resolution $B(A,\,C)$  of a 
 left $A$-module $C$ (see \cite[X.2]{Maclane:63}).
The inverse of the canonical isomorphism $\bar{A} \isor{} A/K\cdot 1_A$, \  $a \mapsto \bar{a}$, yields an isomorphism of   $B(A,\,C)$  with the complex $\T{B}(A,\,C)$ defined by $\T{B}_i(A,\,C) = A\ot_K\bar{A}\htt{i}\ot_K C$ and the same formula for the  differential as in $B(A,\,C)$. But denoting by $\Sigma X$, for a complex $X$,  the suspension of  $X$, we see that $P_nB(A) = \Sigma \T{B}(A/\bar{A}^{n+1},\,\bar{A}/\bar{A}^{n+2} )$, whence
 \[H_i(P_nB(A)) = \left\{ \begin{array}{cl} \ruled\bar{A}/\bar{A}^{n+2}  &\mbox{if $ i=1$}\\ 0 & \mbox{otherwise}    \end{array}\right.\]
and  $P_nH^i(A,\,M)$, $P_nH_i(A,\,N)$ are relative Ext- and Tor-functors, resp. :
\[P_nH^i(A,\,M) = {\rm Ext}_{(A/\bar{A}^{n+1},\,K)}^{i-1}(\bar{A}/\bar{A}^{n+2},\,M), \]
\[P_nH_i(A,N) = {\rm Tor}^{(A/\bar{A}^{n+1},\,K)}_{i-1}(N,\,\bar{A}/\bar{A}^{n+2}) \]
\par\vspace{.25cm}
{Note that if $K=\mathbb Z$ and $M$ is {\em torsion-free}\/ as a \Z-module, then one can replace the groups ${A}/\bar{A}^{n+1}$ and ${A}/\bar{A}^{n+2}$ in  $P_nB(A)$ by their torsion-free quotients; if, in addition, ${A}/\bar{A}^{n+2}$ is a finitely generated \Z-module, we thus obtain a free suspended resolution of $\bar{\tau}(\bar{A}/\bar{A}^{n+2})$, where $\bar{\tau}(X)$ denotes the quotient modulo \Z-torsion of $X$. Thus, in this case, the relative Ext-functors above turn into usual Ext-functors, and we get an isomorphism 
\BE\label{PnH=Ext} P_nH^i(A,\,M) = {\rm Ext}_{\bar{\tau}(A/\bar{A}^{n+1})}^{i-1}(\bar{\tau}(\bar{A}/\bar{A}^{n+2}),\,M) \EE
Similarly, if $N$ is torsion-free as a \Z-module and ${A}/\bar{A}^{n+2}$ is a finitely generated \Z-module, then
\BE\label{PnH=Tor} \bar{\tau}\big(P_nH_i(A,N)\big) = \bar{\tau}\big({\rm Tor}^{\bar{\tau}(A/\bar{A}^{n+1})}_{i-1}(N,\,\bar{\tau}(\bar{A}/\bar{A}^{n+2}))\big) \EE
In particular, the isomorphisms \REF{PnH=Ext} and \REF{PnH=Tor} hold when $A=\Z[G]$ for a finitely generated group $G$, and this is of particular interest if $G$ is finitely generated torsion-free nilpotent (see \cite{Hartl:96a}, \cite{Habil}).}\V

To relate polynomial (co)homology with usual (co)homology, note that
 the canonical isomorphism 
$\T{B}_i(A,\,K) \:\cong\:  \T{B}_{i-1}(A,\,\bar{A})$ actually is an isomorphism of   complexes
$\T{B}(A,\,K)\:\cong\:\Sigma B(A,\,\bar{A})$, since the last term of the differential of $\T{B}_i(A,\,K)$ vanishes. Thus  natural maps
\BE\label{rhodef}\begin{array}{c} \rho_n^*\mapco P_nH^i(A,\,M) \to H^i(A,\,M), \ruled\\
 \rho_{n*}\mapco H_i(A,\,M) \to P_nH_i(A,\,M)
 \end{array}\EE
are induced by the composite map of complexes
\[ \rho_n\mapco B(A,\,K)  \:\cong\: \T{B}(A,\,K)  \:\cong\: \Sigma \T{B}(A,\,\bar{A}) \Sur{\Sigma r_n} \Sigma \T{B}(A/\bar{A}^{n+1},\,\bar{A}/\bar{A}^{n+2} ) = P_nB(A), \]
where $r_n$ is the tensor product of the canonical projections. Explicitly,
  $$\rho_n(a_0\ot \cdots\ot a_i \ot \lambda) = \lambda  \overline{a_0} \ot \overline{a_1 - \epsilon{a_1}} \ot \cdots \ot 
\overline{a_i - \epsilon{a_i}}$$ in degree $i$.
\V

\begin{lem}\label{rhoinjsurj} For $i=2$, the map $\rho_n^*$ (resp. $\rho_{n\,*}$) in \REF{rhodef} is  injective (resp.\ surjective).
\end{lem}

\proof Let $f\in {\rm Hom}_{A/\bar{A}^{n+1}}(P_nB_2(A),\,M)$ be such that $\rho_n^*[f]=0$. Then $r_n^*(f) = \delta_1^*(g)$ with  $g\in {\rm Hom}_A(\bar{B}_0(A,\,\bar{A}),\,M) = {\rm Hom}_A(A\ot_K\bar{A},\,M)$. For $(a_0,\,a_1,\,a_2)\in A \times \bar{A}\times \bar{A}$, we have $f(\overline{a_0} \ot \overline{a_1} \ot \overline{a_2}) = g(a_0a_1 \ot a_2) - g(a_0\ot a_1a_2)$. One has $g(\bar{A}^{n+1} \ot \bar{A}) = \bar{A}^{n+1} g(1\ot \bar{A}) = 0$, and if $a_1\in \bar{A}^{n+1}$, then $g(a_0\ot a_1a_2)=g(a_0a_1\ot a_2) - f(\overline{a_0} \ot 0\ot \overline{a_2}) =0$. Thus $g$ factors through $\bar{g}\mapco 
A/\bar{A}^{n+1} \ot_K \bar{A}/\bar{A}^{n+2} \to M$, whence $f=\bar{\delta}_2^*(\bar{g}) \equiv 0$ in $P_nH^2(A,\,M)$ and $\rho_n^*$ is injective. As to $\rho_{n\,*}$, note that for $(x,\,a_1,\,a_2)\in N\times \bar{A}^{n+1} \times \bar{A}$, we have 
$(N\ot \delta_2)(x\ot a_1\ot a_2) = xa_1\ot a_2 - x\ot a_1a_2 = - x\ot a_1a_2$
since $N$ is an $A/ \bar{A}^{n+1}$-module. Thus $\Imm{N\ot  \bar{A}^{n+2}} = (N\ot \delta_2)\Imm{N\ot 
\bar{A}^{n+1}\ot \bar{A}}$ and whence $\Ker{
N\ot \bar{\delta}_2} = (N\ot \rho_2)\Ker{N\ot \delta_2}$. This implies that $\rho_{n\,*}$ is surjective.
\hfbox\V

The $K$-modules $P_nH^i(A,\,M)$ and $P_nH_i(A,\,N)$, for fixed $i\ge 2$ and varying $n\ge 0$, are related by   chains of natural maps
\[ 0 =P_0H^i(A,\,M) \to   \ldots \to P_nH^i(A,\,M) \mr{\sigma_n^*} P_{n+1}H^i(A,\,M) \to \ldots \mr{\rho^*} H^i(A,\,M)\]
\BE\label{H2cofiltr} H_i(A,\,N) \mr{\rho_*} \ldots  \to P_{n+1}H_i(A,\,N) \mr{\sigma_{n*}} P_{n}H_i(A,\,N) \to \ldots \to P_{0}H_i(A,\,N) =0\ruled \EE
commuting with the maps $\rho^*$ and $\rho_*$ where $\sigma_n\mapco P_{n+1}B(A) \auf P_{n}B(A)$ is the tensor product of the canonical projections.
For $i=2$ the map $\sigma_n^*$  (resp. $\sigma_{n*}$) is injective (resp.\ surjective) by  Lemma \ref{rhoinjsurj}; therefore,  identifying $P_nH^2(A,\,M)$ with its isomorphic image $\rho_n^*P_nH^2(A,\,M)$ in $H^2(A,\,M)$, provides a {\em natural ascending filtration}\/ of $H^2(A,\,M)$,
\[ 0=P_0H^i(A,\,M) \subset   \ldots  \subset P_nH^i(A,\,M)  \subset P_{n+1}H^i(A,\,M)  \subset \ldots  \subset H^i(A,\,M)\]
Dually, the maps $\sigma_{n*}$ in \REF{H2cofiltr} being surjective for $i=2$, can be interpreted as
  a {\em natural cofiltration}\/ of $H_2(A,\,N)$ which, in turn, gives rise to a {\em natural descending filtration}\/  
\BE\label{Kerrhofiltr}H_2(A,\,N) =\Ker{\rho_{0*}}  \supset \ldots \supset \Ker{\rho_{n*}} \supset \Ker{\rho_{n+1*}}   \supset \ldots \supset 0\EE
Now, returning to polynomial (co)homology of groups, note that $$P_nB_i(\Z[G]) = \Z[G]/\mathfrak g^{n+1} \ot P_n(G)\htt{i-1}\ot P_{n+1}(G),$$ whence 
$\rho_n^*P_nH^i(G,\,M)$ is the subgroup of $H^i(G,\,M)$ consisting of elements representable by multi-polynomial cocycles of degree $\le n$ in the first  $i-1$ variables and of degree $\le n+1$ in the last variable. For $i=2$, note that 
 \begin{eqnarray}\label{coeffZ}
& & \hspace{-10mm}\coker{\Z\ot \bar{\delta}_3\mapco \Z\ot_G P_nB_3(\Z(G)) \to \Z\ot_G P_nB_2(\Z(G))} \nonumber\\ &\cong& P_n(G) \ot_G P_{n+1}(G) \nonumber\\
&\cong& P_n(G) \ot_G P_{n}(G)
\end{eqnarray}
 Thus, if $M$ is a {\em trivial }\/ $G$-module, 
$\rho_n^*P_nH^2(G,\,M)$ is the subgroup of $H^2(G,\,M)$ consisting of elements representable by bipolynomial cocycles of degree $\le n$ in both variables; consequently, $P_nH^2(G,\,M)$ is isomorphic via $\rho_n^*$ with the polynomial cohomology groups defined in \cite{Passi:68b}.
\para
We now wish to show that  actually it is the filtration \REF{Kerrhofiltr}
which is intimately related to dimension subgroups. Denote by 
\[ \rho_{n*}^G = \rho_{n*} \mapco H_2(G) \to P_nH_2(G), \]
abbreviating $H_2(G)=H_2(G,\,\Z)$ and $P_nH_2(G)
=P_nH_2(G,\,\Z)$.\V

\begin{lem}\label{rhoexpl} One has a commutative diagram with exact rows
\[\begin{matrix}
\ruled 0&\to &H_2(G) &\mr{j}& \mathfrak g\ot_G \mathfrak g & \mr{\mu} &\mathfrak g \cr
& & \mapdown{\rho_{n*}^G} & & \mapdown{\pi_n \ot \pi_n}& & \mapdown{\pi_{n+1}} \cr 
\ruleu 0&\to &P_nH_2(G) &\mr{j_n}&  P_n(G) \ot_G P_n(G) & 
\mr{\mu_n} & P_{n+1}(G),
\end{matrix}\]
where $\mu,\,\mu_n$ are induced by multiplication in $\mathfrak g$ and $\pi_k$ is the canonical projection.
\end{lem}
\par\vspace{.25cm}
The above lemma easily follows from \REF{coeffZ}. We need to specify an isomorphism in the Hopf formula. We use the commutator conventions $[a,\,b]=aba^{-1}b^{-1}$ for $a,\,b\in G$ and $[x,\,y]=xy-yx$ for elements $x,\,y$ in some ring. As usual, $G'=\gamma_2(G)$. For the rest of this section, choose a free presentation
 \[   1\to R\mr{i'} F \mr{q'} G \to 1 \]
 of $G$, where we consider $R$ as a subgroup of $F$. It is often convenient to write $\bar{a} =q'(a)$ for $a\in F$.\V

\begin{prop}\label{Hopf} An isomorphism 
\[ \nu \mapco R\cap F'/[F,\,R] \hspace{2mm}\isor{}\hspace{2mm} H_2(G) = \Ker{\mu \mapco \mathfrak g\ot_G \mathfrak g \to \mathfrak g} \]
is given by the composite map $R\cap F'/[F,\,R] \hra [F,\,F]/[F,\,R] \mr{\nu'} \mathfrak g\ot_G \mathfrak g$, where for $a,\,b\in F$  
\begin{multline}
\nu'(\overline{[a,\,b]}) =  (\bar{a}-1) \ot (\bar{b}-1)\bar{a}^{-1}\bar{b}^{-1} -(\bar{b}-1) \ot (\bar{a}-1)\bar{a}^{-1}\bar{b}^{-1} \label{nudefvar}\\
= (\bar{a}-1) \ot (\bar{b}-1)  -(\bar{b}-1) \ot (\bar{a}-1) + [\bar{a}-1\,,\,\bar{b}-1]\ot  (\bar{a}^{-1}\bar{b}^{-1} - 1).   \end{multline}
\end{prop}
\par\vspace{.25cm}\noindent
\proof 
Consider the following  diagram of homomorphisms.
\BE\label{kappadia1}\begin{matrix}
\ruleu\ruled 1 & \to & R\cap F'/[F,\,R] & \mr{\nu} & \mathfrak g\ot_G \mathfrak g & \mr{\mu} & \mathfrak g^2 & \to & 0\cr
& & \injdown{\iota} & &\mapdown{\mu'} & & \injdown{\iota} \cr
\ruled1 & \to & R/[F,\,R] & \mr{D\overline{i'}} & I\Big(F/[F,\,R] \,,\,R/[F,\,R] \Big) & \mr{I(\overline{q'})} & \mathfrak g & \to & 0\cr
& & \surdown{q} & & \surdown{\chi} & &\surdown{\chi} \cr
\ruleu1 &\to& R/R\cap F' &\hra & F/F'& \mr{\overline{q'}} &G/G' &\to& 1\ruled
\end{matrix}\EE
Here $\iota$ denotes the respective inclusions, $q$ the natural projection, the maps $\chi$ send $\overline{a-1}$ to $a\Gamma'$ for $a\in \Gamma=F,\,G$, and
$\mu'( (\bar{a}-1)\ot (\bar{b}-1)) = \overline{(a-1)(b-1)}$, $a,\,b\in F$, is well-defined by centrality of $R/[F,\,R]$.
The diagram commutes; for the upper left hand square this follows from the identity $[a,\,b]-1 = [a-1,\,b-1]a^{-1}b^{-1}$ for $a,\,b\in F$.
The middle and bottom row are exact, and so are the columns, thanks to the well-known isomorphism $\chi\mapco P_1(\Gamma)\hcong\Gamma/\Gamma'$ for any group $\Gamma$. This successively implies  injectivity of
 $\mu'$, the relation
   \BE\label{cap}  (D\overline{i'})^{-1} \Imm{\mu'} = R\cap F'/[F,\,R] \EE
and exactness of the first row, as asserted.\hfbox\V

Now we are ready for the main result of this section.
\par\vspace{.25cm}
\begin{satz}\label{DnH2} Let $e\mapco  C\stackrel{i}{\hra} E \stackrel{q}{\auf} G$ be a central group extension and suppose that $\gamma_n(E)=1$. Then 
\[ D_n(E,\,C) \cap C = \kappa\Big(\Ker{\rho^G_{n-2\,*}\mapco H_2(G) \auf P_{n-2}H_2(G)} \Big)\]
where $\kappa\mapco H_2(G) \to C$ is adjoint to the cohomology class of $e$ under the Kronecker pairing $H^2(G,\,C) \times H_2(G) \to C$.
\end{satz}

\par\vspace{.25cm}\proof This is just an extension of the proof of Proposition \ref{Hopf}. With the notations there 
we obtain  a commutative diagram of central extensions
\[\begin{matrix}
\ruled 1 & \to & R/[F,\,R] & \mr{\overline{i'}} & F/[F,\,R] & \mr{\overline{q'}} & G & \to & 1 \cr
& & \mapdown{\alpha_1} & & \mapdown{\alpha_0} & & \Big{\|} \cr
\ruleu 1 & \to & C & \mr{ {i}} & E & \mr{{q }} & G & \to & 1\ruled
\end{matrix}\]
by choosing any lifting  $\alpha_0$ of $\overline{q'}$.  Then $\kappa =\alpha_1\iota \nu^{-1}$. Consider the following commutative diagram with exact rows.
\BE\label{kappadia}\begin{matrix}
\ruled 0 & \to & H_2(G) & \mr{j} & \mathfrak g\ot_G \mathfrak g & \mr{\mu} & \mathfrak g^2 & \to  & 0 \cr
& & \mapdown{\iota\nu^{-1}} & &\mapdown{\mu'} & &\mapdown{\iota}\cr
\ruled1 & \to & R/[F,\,R] & \mr{D\overline{i'}} & I\Big(F/[F,\,R] \,,\,R/[F,\,R] \Big) & \mr{I(\overline{q'})} & \mathfrak g & \to & 0\cr
& & \mapdown{\alpha_1} & & \mapdown{I(\alpha_0)} & &\Big{\|} \cr
\ruleu1 &\to& C &\mr{Di} &I(E,\,C)& \mr{I(q)} &\mathfrak g &\to& 0\ruled
\end{matrix}\EE
Note that $x\in H_2(G)$ lies in \Ker{\rho_{n-2\,*}^G} if and only if  $j(x) \in \Imm{\mathfrak g\ot \mathfrak g^{n-1}}$. Thus $Di\kappa(\Ker{\rho_{n-2\,*}^G}) \subset P_n(E,\,C)$ and $\kappa(\Ker{\rho_{n-2\,*}^G}) \subset D_n(E,\,C)$. 
\para
Conversely, let $c\in D_n(E,\,C) \cap C$. Then there is $y\in  \Imm{\mathfrak g\ot \mathfrak g^{n-1}}$ such that $Di(c) = I(\alpha_0)\mu'(y)$. By exactness of the rows of the diagram \REF{kappadia} it follows that there is $z\in R/[F,\,R]$ such that $D\overline{i'}(z) =\mu'(y)$. By \REF{cap} there exists  $x\in H_2(G)$ such that $\iota\nu^{-1}(x) = z$. Then $\mu'j(x) = \mu'(y)$, whence $j(x)=y$, since $\mu'$ is injective. Hence $x\in \Ker{\rho_{n-2\,*}^G}$. But $Di\kappa(x) = I(\alpha_0)\mu'(y) = Di(c)$, whence $c=\kappa(x) \in \kappa(\Ker{\rho_{n-2\,*}^G})$, as desired.\hfbox\V

\begin{kor}\label{redHniln-2} Let $E$ be an $(n-1)$-step nilpotent group and $C$  a central subgroup of $E$. Then $D_n(E,\,C) \cap C\gamma_{n-1}(E)$ is a homomorphic image of \Ker{\rho_{n-2\,*}^G} for the $(n-2)$-step nilpotent group $G=E/C\gamma_{n-1}(E)$.
\end{kor}
\par\vspace{.25cm}
In fact,  $D_n(E,\,C) =  D_n(E,\,C\gamma_{n-1}(E)) $, since $\gamma_{n-1}(E) -1 \subset \mathfrak e^{n-1}$.\V

As a first  illustration of the method, we reprove a result known for a long time (see \cite{Passi:79}).\V

\begin{kor}\label{D2D3} For any group $\Gamma$ and $n\le 3$, $D_n(\Gamma)=\gamma_n(\Gamma)$.
\end{kor}

\par\vspace{.25cm}\proof We apply Corollary \ref{redHniln-2} to $E=\Gamma/\gamma_n(\Gamma)$. The case $n=1$ is trivial.  For $n=2$, take   $C=E$. Then $D_2(\Gamma)/\gamma_2(\Gamma) = D_2(E,\,E)\cap E$ is trivial as $G$ is. 
\para
For $n=3$, take $C=\gamma_2(E)$. Then $$D_3(\Gamma)/\gamma_3(\Gamma) =D_3(E)= D_3(E,\,\gamma_2(E)) \cap \gamma_2(E),$$ since $\gamma_2(E)-1\subset \mathfrak e^2$ and $D_3(E) \subset D_2(E) =\gamma_2(E)$. 
Consider the  maps $$G\sm G \mr{\phi} R\cap F'/[F,\,R]=F'/[F,\,R] \mr{\nu} \mathfrak g\ot_G \mathfrak g \Sur{\pi_1\ot\pi_1} P_1(G) \ot_G P_1(G) = G \ot G,$$ where $\phi$ is the classical isomorphism given by $j\phi (\bar{a}\sm \bar{b}) = \overline{[a,\,b]}$ for $a,\,b\in F$. The composite map $(\pi_1\ot\pi_1)\nu\phi$ sends $a\sm b$ to $a\ot b -b\ot a$ by \REF{nudefvar}; but this map is well-known to be injective so $\Ker{\rho_{1\,*}^G}=0$ by Lemma \ref{rhoexpl}.\hfbox\V

\section{The fourth dimension quotient}\label{fourthdq1}
\para
We now consider the case $n=4$ in Corollary \ref{redHniln-2}. This leads to a  simple description of 
the crucial group $\Ker{\rho_{2\,*}^G}$ in terms of the exterior torsion square of $G\ab=G/G'$ which amounts to a nice proof of the well-known  result that, for all groups $\Gamma$, the dimension quotient $D_4(\Gamma)/\gamma_4(\Gamma)$ is of exponent $2$ (see \cite{Passi:79}).\V


For an Abelian group $A$, let $L(A),\,T(A)$, and $S(A)$ denote the free Lie algebra, the tensor algebra and the symmetric algebra over $A$  respectively. The natural maps of graded Abelian groups $L(A) \Inj{l} T(A) \Sur{s} S(A)$ are the injection into the universal enveloping algebra and the canonical projection, resp.  Thus $l_n$ sends an $n$-fold Lie bracket in $L_n(A)$ to the corresponding tensor commutator in $T_n(A)=A\htt{n}$. In particular, $L_2(A)=A\sm A$ and $l_2(a\sm b) = a\ot b -b\ot a$ for $a,\,b\in A$.
\para
For the rest of this section, let $G$ be a {\em $2$-step nilpotent group}\/.   The surjective homomorphism $c_2\mapco G\ab \sm G\ab \auf G'$ is defined by $c_2(\bar{a} \sm \bar{b}) = [a,\,b]$ for $a,\,b\in G$. For $x\in G\ab$ and $m\in \Z$ such that $mx=0$, choose elements $\T{x}\in G$ and $f_mx \in G\ab\sm G\ab$ such that $\T{x}G'=x$ and $c_2(f_mx) = \T{x}^m$.\V

The main ingredients for calculating $\Ker{\rho_{2\,*}^G}$ are the structure theorems describing $H_2(G)$ and $P_2(G)\ot_GP_2(G)$ for 2-step nilpotent  $G$ which we
recall  from \cite{Hartl:96} and \cite{Hartl:95}.\V

\begin{satz}\label{H2PoP} For 2-step nilpotent $G$ there are natural exact sequences
\BE\label{H2sequ}
{\rm Tor}(G\ab,\,G\ab) \hmr{\delta}  \frac{\dst L_3(G\ab)}{\dst [G\ab,\,\Ker{c_2}] + V} \hmr{\nu i} H_2(G) \hmr{\sigma} G\ab \sm G\ab \hmr{c_2} G' \to 1,
\EE
\[
{\rm Tor}(G\ab,\,G\ab) \hmr{\delta'}  \frac{\dst (G\ab)\htt{3}}{\dst l_2\Ker{c_2}\ot G\ab +G\ab\ot l_2\Ker{c_2}} \hmr{i'} P_2(G)\ot_G P_2(G)\] 
\BE\label{PoPsequ} \hspace{9.5cm} \hmr{\sigma'}G\ab\ot G\ab \to 0,
\EE
where \[\delta (\tau_m(x_1,\,x_2)) = q\Big( [x_1,\,  f_m x_2]  + [x_2,\,  f_m x_1]   + \binom{m}{2} [x_1+x_2,\,[x_1,\,x_2]] \Big)\]
\[\delta' (\tau_m(x_1,\,x_2)) = q'\Big( x_1\ot (l_2f_m x_2) - 
(l_2f_m x_1) \ot x_2 + \binom{m}{2} (x_1\ot x_1 \ot x_2 -
x_1\ot x_2 \ot x_2) \Big)\]
with $q,\,q'$ being the canonical projections,
and $$i[\bar{a},\,[\bar{b},\,\bar{c}]] = \overline{[a,\,[b,\,c]]},$$ 
$$\ruled \sigma \mapco H_2(G) \isor{\nu^{-1}} \frac{R\cap F'}{[F,\,R]} \hra \frac{[F,\,F]}{[F,\,R]} \mr{\nu''} G\ab\sm G\ab$$ with $$\nu''(\overline{[a,\,b]}) = \bar{a} \sm \bar{b},$$ $$i'(\bar{a}\ot \bar{b}\ot \bar{c})= \overline{(a-1)(b-1)} \ot \overline{(c-1)},$$ $$\sigma'(\overline{(a-1)} \ot \overline{(b-1)}) = \bar{a}\ot \bar{b},$$ for $a,\,b,\,c\in G$. Finally, $V$ denotes the subgroup of $L_3(G\ab)$ generated by the elements $[x,\,f_{o(x)}x]$, where $x$ ranges over the elements of finite even order $o(x)$ of $G\ab$.\hfbox

 \end{satz}
\par\vspace{.25cm}
Note that for any torsion element $x$ of $G\ab$, $\delta(\tau_{o(x)}(x,\,x)) = 2 [x,\,f_{o(x)}x]$, so if $o(x)$ is odd, $[x,\,f_{o(x)}x] \in \Imm{\delta}$. Thus $V$ can be replaced by the subgroup $V'$ generated by the elements $[x,\,f_{o(x)}x]$ for {\em any }\/ torsion elements $x\in G\ab$. Now if $mx=0$ for $m\in \Z$, then $[x,\,f_m x]\in V'$; this shows that the map
\[ \delta_1\mapco G\ab \stackrel{\sm}{*}G\ab \hmr{} 
\frac{\dst L_3(G\ab)}{\dst [G\ab,\,\Ker{c_2}] + V + \Imm{\delta}}\]
defined by $\delta_1(\overline{  \tau_m(x_1,\,x_2) }) = \overline{[x_2,\,f_mx_1]}$ is well-defined. Moreover, define homomorphisms
\[  G\ab \ot G' \hml{\delta_2}  G\ab \stackrel{\sm}{*}G\ab \hmr{\delta_3} SP^3(G\ab) \]
by  $\delta_2(\,\overline{  \tau_m(x_1,\,x_2) }\,) = x_1\ot \T{x}_2^k - x_2 \ot \T{x}_1^k$ and $\delta_3(\,\overline{ \tau_m(x_1,\,x_2)}\,) = \binom{m}{2} s_3(x_1\ot x_1 \ot x_2 - x_1\ot x_2 \ot x_2)$. 
\para
Our  main result of this section is the following: 
\para

\begin{satz}\label{Kerrho2} For every  2-step nilpotent group $G$,  
\[  \Ker{\rho_{2\,*}^G} = \overline{\nu i} \delta_1 (\Ker{\delta_2} \cap \Ker{\delta_3} ) \]
\end{satz}
\par\vspace{.25cm}
We note that this result, combined with Corollary \ref{redHniln-2},
 leads to the construction of a multi-parameter family of examples of groups with non-trivial fourth dimension quotient \cite{Hartl:D4}.
\para
Before proving the theorem, we note an immediate consequence which in view of Corollaries \ref{redHniln-2} and \ref{D2D3}  reproves the known result stating that $D_4(E,\,C)/\gamma_4(E)$ is of exponent $1$ or $2$.
\par\vspace{.25cm}
\begin{kor}\label{expo2} For every  2-step nilpotent  group $G$,  $$2 \Ker{\rho_{2\,*}^G} =0\ \text{and}\ \Ker{\rho_{2\,*}^G} \subset 2 \Imm{\overline{\nu i}} \subset 2 H_2(G).$$
\end{kor}

\par\vspace{.25cm}\proof Let $\beta\mapco G\ab \ot G' \to\frac{\dst L_3(G\ab)}{\dst [G\ab,\,\Ker{c_2}] + V + \Imm{\delta}}$ be the factorisation of the bracket map $G\ab \ot L_2(G\ab) \to L_3(G\ab)$ followed by the projection through $1\ot c_2\mapco G\ab \ot L_2(G\ab) \auf G\ab\ot G'$. Now
\begin{eqnarray*}
-\beta \delta_2 (\overline{ \tau_m(x_1,\,x_2)}) &=& \overline{{}-[x_1,\,f_m x_2] + [x_2,\,f_m x_1]} \ruled\\
&=& \delta_1(\overline{\tau_m(x_1,\,x_2)} - \overline{ \tau_m(x_2,\,x_1)}) \ruled\\
&=& \delta_1(2 \overline{ \tau_m(x_1,\,x_2)})
\end{eqnarray*}
whence $2\delta_1 = -\beta\delta_2$ and $2\delta_1(
\Ker{\delta_2} \cap \Ker{\delta_3} ) =0$. The remaining relations are immediate from Theorem \ref{Kerrho2} and Lemma \ref{Kerdel3} below.\hfbox\V

\begin{lem}\label{Kerdel3} The subgroup \Ker{\delta_3} of $G\ab \stackrel{\sm}{*}G\ab$ is generated by the elements $\overline{\tau_m(x_1,\,2x_2)}$,\, $x_1,\,x_2\in G\ab$, $m\in \Z$ such that $mx_1=2mx_2=0$.
\end{lem} 

\par\vspace{.25cm}\proof It is clear that the indicated elements are in \Ker{\delta_3}. To prove the converse, we may assume that $G\ab$ is a finite $2$-group. Let $G\ab = \bigoplus_{i=1}^n\Z/2^{r_i} \cdot x_i$, $1\le r_i\le r_j$ if $i\le j$, be a decomposition of $G\ab$ into  cyclic factors generated by elements $x_i$. Then $G\ab \stackrel{\sm}{*}G\ab = \bigoplus_{1\le i<j\le n} \Z/2^{r_i} \overline{\tau_{2^{r_i}} (x_i,\, 2^{r_j-r_i} x_j )}$. One has 
\begin{eqnarray*}
\delta_3 \overline{\tau_{2^{r_i}} (x_i,\, 2^{r_j-r_i} x_j )} &=& 
2^{r_i-1}(2^{r_i} -1) s_3\big(x_i\ot x_i\ot 2^{r_j-r_i}x_j - x_i
\ot 2^{r_j-r_i}x_j \ot 2^{r_j-r_i}x_j \big) \\
&=&   2^{r_j-1} s_3(x_i\ot x_i \ot x_j) - 2^{2r_j-r_i-1}s_3(x_i\ot x_j \ot x_j )
\end{eqnarray*}
So $\delta_3 \overline{\tau_{2^{r_i}} (x_i,\, 2^{r_j-r_i} x_j )} =0$ if $r_i<r_j$. As $SP^3(G\ab) = \bigoplus_{1\le i\le j\le k\le n} \Z/2^{r_i} \cdot s_3(x_i \ot x_j \ot x_k) $ we see that the above cyclic factors of 
$G\ab \stackrel{\sm}{*}G\ab$ map into different direct components of $SP^3(G\ab)$ under $\delta_3$; hence \Ker{\delta_3} is generated by the elements $
\overline{\tau_{2^{r_i}} (x_i,\, 2^{r_j-r_i} x_j )}$, $r_i<r_j$, and $2\overline{\tau_{2^{r_i}}( x_i,\, x_j)} = \overline{\tau_{2^{r_i}}( x_i,\, 2x_j)}$, $r_i=r_j$.
 \hfbox\V
\newpage
\proofof{Theorem \ref{Kerrho2}} Consider the following diagram.
\BE\label{22dia} 
\begin{matrix}
{\rm Tor}(G\ab,\,G\ab)  \oplus V''  &  \mr{\delta_V}  & \frac{\ruleu\dst L_3(G\ab)}{\ruled\dst K}  &  \mr{\nu i}  &  H_2(G)  &  \mr{\sigma }  & \Ker{c_2}  &  \to & 0\cr
\mapdown{t}  &  &  \mapdown{\overline{l_3}} & & \mapdown{\rho^G_{2*}}  & & \injdown{l_2} & & \cr
{\rm Tor}(G\ab,\,G\ab)  & \mr{\delta'}  &  
\frac{\ruleu\dst (G\ab)\htt{3}}{\ruled\dst K'}  &\hmr{i'}&  P_2(G)\ot_G P_2(G)  & \mr{\sigma'}  & G\ab\ot G\ab  & \to & 0\cr
\surdown{\pi}  &  & \surdown{\pi'}  &  & \|  & & \|  &  & \cr
\ruleu\ruled G\ab \stackrel{\sm}{*}G\ab &  \mr{\overline{\delta'}}  &  \coker{\overline{l_3}\delta_V}  &  \mr{\overline{i'}}  & P_2(G)\ot_G P_2(G)  & \mr{\sigma'}  & G\ab\ot G\ab  & \to & 0\cr
\injup{j}  &  &  \injup{\overline{\overline{l_3}}}  &  &  \mapup{\rho^G_{2*}}  &  & \injup{l_2}  &   &  \cr
\Ker{\delta_2} \cap \Ker{\delta_3}    & \mr{\delta_1}  & \frac{\ruleu\dst L_3(G\ab)}{\ruled\dst K+\Imm{\delta_V}}  &  \Inj{\overline{\nu i}}  & 
H_2(G)  &  \mr{\sigma }  & \Ker{c_2}  &  \to & 0

\end{matrix}\EE

\N Here $V''$ is the free Abelian group generated by elements $[x]$ where $x$ runs through the torsion elements of $G\ab$. Moreover, $K= [G\ab,\,\Ker{c_2}]$ and $K'=l_2\Ker{c_2}\ot G\ab +G\ab\ot l_2\Ker{c_2}$. The map $\delta_V$ is $\delta$ on ${\rm Tor}(G\ab,\,G\ab)$ and sends $[x]$ to $[x,\,f_{o(x)}x]$, while $t$ sends $\tau_m(x_1,\,x_2)$ to 
$\tau_m(x_1,\,x_2) + \tau_m(x_2,\,x_1)$ and $[x]$ to $\tau_{o(x)}(x,\,x)$.
The maps $\pi,\,\pi'$ are the canonical projections while $j$ is the canonical injection.

It is straightforward to check that the three upper squares commute; 
and the two upper rows of diagram \REF{22dia} are exact by Theorem \ref{H2PoP}.  This implies that $\overline{\delta'}$ and $\overline{i'}$ are well-defined, that the third row is exact, too, and that all squares commute, except possibly  the one in the bottom left-hand corner. By injectivity of $l_2$ we have
\BE\label{Kerrhoids}
\Ker{\rho^G_{2*}} =  \overline{\nu i}\,\overline{\overline{l_3}}^{\,-1} (\Imm{\overline{\delta'}} \cap \Imm{ \overline{\overline{l_3}}} )
\EE
Now we recall the well-known natural exact sequence for any Abelian group $A$ (which is a natural version of the Poincar\'e-Birckhoff-Witt decomposition of $A\htt{3}$)
\[ 0 \mr{} L_3(A) \oplus A \ot (A\sm A) \mr{(l_3,\,1\ot l_2)} A\htt{3} \mr{s_3} SP^3(A) \mr{} 0 \]
which for $A=G\ab$ successively induces the following exact sequences natural in $G$
\[ 
0 \mr{} \frac{\dst L_3(A)}{\dst [G\ab,\,\Ker{c_2}]} \hspace{1mm}\oplus\hspace{1mm}  \Big( G\ab \ot \frac{\dst (G\ab\sm G\ab)}{\dst \Ker{c_2}} \Big)  \mr{(\overline{l_3},\,\overline{1\ot l_2)}} \frac{\dst (G\ab)\htt{3}}{K'} \mr{\overline{s_3}} SP^3(A) \mr{} 0 \]
\BE \label{PBW}
0 \to \frac{\dst L_3(A)}{\dst [G\ab,\,\Ker{c_2}]+\Imm{\delta_V}} \hspace{2mm}\oplus\hspace{2mm}   ( G\ab \,\ot\, G')  \mr{(\overline{\overline{l_3}},\,\pi\overline{1\ot l_2)}\, } \coker{\overline{l_3}\, \delta_V} \mr{\overline{\overline{s_3}}} SP^3(A) \to 0 \EE
One has $\Imm{\overline{\delta'}} \cap \Imm{ \overline{\overline{l_3}}} \subset \overline{\delta'}\Ker{\overline{\overline{s_3}} \,\overline{\delta'}}$. But $\overline{\overline{s_3}} \,\overline{\delta'} =\delta_3$, and on \Ker{\delta_3}, $ \overline{\delta'}$ coincides with the map $\delta''$ sending $\overline{\tau_m(x_1,\,x_2)}$ to $\pi q'([x_2,\,l_2f_m x_1] + x_1 \ot l_2f_m x_2 - x_2 \ot l_2f_m x_1 ) = ( \overline{\overline{l_3}} \delta_1 +  \pi\overline{1\ot l_2} \delta_2 ) (\,\overline{\tau_m(x_1,\,x_2)}\,)$ since we can write $\delta'$ as $\delta'(\tau_m(x_1,\,x_2)) = 
[x_2,\,l_2f_m x_1] + x_1 \ot l_2f_m x_2 - x_2 \ot l_2f_m x_1 + \binom{m}{2}(x_1\ot x_1 \ot x_2 -
x_1\ot x_2 \ot x_2)$, and  by Lemma \ref{Kerdel3}.
Thus the square in  bottom left-hand corner of diagram \REF{22dia} commutes.
By sequence \REF{PBW}, $\Imm{\overline{\delta'}} \cap \Imm{ \overline{\overline{l_3}}} = 
\overline{\overline{l_3}}\delta_1(\Ker{\delta_3}\cap \Ker{\delta_2})$, whence 
$\Ker{\rho^G_{2*}} =  \overline{\nu i}\,\delta_1
(\Ker{\delta_3}\cap \Ker{\delta_2})$ by injectivity of $
\overline{\overline{l_3}}$, as asserted.\hfbox
\Para
\section{Two spectral sequences}\label{ss}

\Para Let {\sf Gr} denote the category of groups. Using simplicial
homotopy theory, for every functor $T:{\sf Gr}\to {\sf Gr}$, one
can define its left derived functors $\mathcal L_qT:{\sf Gr}\to
{\sf Gr}$ (see \cite{Keune:73}). If $\mathbf F\to G$ is a free
simplicial resolution of $G$, then $\mathcal
L_qT(G)=\pi_q(T(\mathbf F))$, the $q$-th homotopy group of the
simplicial group $T(\mathbf F)$.

\para
For $G\in {\sf Gr}$, as usual, let $\{\gamma_n(G)\}_{n\ge 1}$ denote its lower central series. We then have functors $$\varGamma_n:{\sf Gr}\to {\sf Gr},\quad \mathfrak L_n:{\sf Gr}\to {\sf Gr},\quad n\ge 1,$$given by $\varGamma_n(G)=G/\gamma_n(G),\ \mathfrak L_n(G)=\gamma_n(G)/\gamma_{n+1}(G)$, which extend naturally to the category ${\sf sGr}$ of simplicial groups.  In view of the natural exact sequences $$1\to \mathfrak L_n(G) \to \varGamma_{n+1}(G)\to \varGamma_n(G)\to 1,\quad n\ge 1,\quad G\in {\sf Gr},$$
we have a homotopy exact couple associated with the exact sequences
$$1\to \mathfrak L_n(\mathbf F)\to \varGamma_{n+1}(\mathbf F)\to \varGamma_n(\mathbf F)\to 1, \quad n\geq 1,$$of simplicial groups. We thus have a spectral sequence $E(G):=\{E^r_{p,\,q}(G)\}$, with $$E^1_{p,\,q}(G)=\mathcal L_q\mathfrak L_p(G),$$ and the differential $d^r$ having degree $(r,\,-1)$.
\par\vspace{.5cm}
Let $\mathfrak P_n:{\sf Gr}\to {\sf Gr}$ and $\mathfrak  Q_n:{\sf Gr}\to {\sf Gr}$ be the functors defined  by $$\mathfrak P_n(G)=\mathbb Z[G]/\mathfrak g^n,\quad \mathfrak Q_n(G)=\mathfrak g^n/\mathfrak g^{n+1},\quad n\ge 1,$$ We then have natural exact sequences $$0\to \mathfrak Q_n(G)\to \mathfrak P_{n+1}(G)\to \mathfrak P_n(G)\to 0,\quad n\ge 1,$$ The homotopy exact couple associated with the induced  exact sequences
$$0\to \mathfrak Q_n(\mathbf F)\to \mathfrak P_{n+1}(\mathbf F)\to\mathfrak  P_n(\mathbf F)\to 0,\quad n\geq 1,$$
 yield another spectral sequence $\overline{E}(G):=\overline
E_{p,\,q}^r(G)$ with $$\overline E_{p,\,q}^1(G)=\mathcal
L_q\mathfrak Q_p(G),$$ and the differential $\overline d^r$  again having
bidegree $(r,\,-1)$. There exists  a natural homomorphism
$$\kappa:E(G)\to \overline{E}(G)$$induced by the canonical
injection $\kappa:G\to \mathbb Z[G],\ g\mapsto g-1,\ g\in G$. It
is known that $$ E_{n,\,0}^\infty(G)=\gamma_n(G)/\gamma_{n+1}(G),\
\ \overline E_{n,\,0}^\infty(G)=\mathfrak g^n/\mathfrak g^{n+1},\
n\geq 1.
$$
For more details about these two spectral sequences, see \cite{Gruenenfelder:80}.

\para It is clear (see \cite{Kan:58}) that
$$
E_{1,\,q}^1(G)=\overline E_{1,\,q}^1(G)=H_{q+1}(G),\ q\geq 0.
$$

The above spectral sequences define certain filtrations in group
homology:
\begin{align*}
& H_{m+1}(G)=E_{1,\,m}^1(G)\supseteq E_{1,\,m}^2(G)\supseteq
E_{1,\,m}^3\supseteq \dots\\
& H_{m+1}(G)=\overline E_{1,\,m}^1(G)\supseteq \overline
E_{1,\,m}^2(G)\supseteq \overline E_{1,\,m}^3\supseteq \dots
\end{align*}
It may be observed that, for the case $m=1$, these filtrations are
the Dwyer's filtration \cite{Dwyer:75} and the dual of the
filtration studied by Passi-Stammbach \cite{PS:74} respectively, the latter being the same as the filtration (\ref{Kerrhofiltr}).

\Para {\bf The initial terms.} Let $1\to R\to F\to G\to 1$ be a
free presentation of the group $G$, and  $\mathbf F\to G$  a
free simplicial resolution of $G$ with $F_0=F$.
\par\vspace{.25cm}
We make the following standard notations:
\begin{align*}
& R(0)=R,\\
& R(k+1)=[R(k),F],\ k\geq 0,\\
& {\bf r}(0)=(R-1)\mathbb Z[F],\\
& {\bf r}(k+1)=\mathfrak f{\bf r}(k)+{\bf r}(k)\mathfrak f,\ k\geq
0.
\end{align*}
Observe that for every $k\geq 1$ and $x\in R(k),$ one has $1-x\in
{\bf r}(k)$. Since Lie functors and universal enveloping functors
preserve coequalizers, direct simplicial computations imply that
\begin{align*}
& E_{n,\,0}^1(G)=\gamma_n(F)/R(n-1)\gamma_{n+1}(F),\ n\geq 1,\\
& \overline E_{n,\,0}^1(G)=\mathfrak f^n/({\bf r}(n-1)+\mathfrak
f^{n+1}),\ n\geq 1,
\end{align*}
and the map $\kappa$ is a monomorphism on the lower level
\cite{Gruenenfelder:80}:
\begin{equation}\label{Gru1}
\kappa_{n,\,0}^1: E_{n,\,0}^1(G)\hookrightarrow \overline
E_{n,0}^1(G),\ n\geq 1.
\end{equation}
This implies the following fact which is due to Sj\"{o}gren
\cite{Sjogren:79}:
\par\vspace{.25cm}
\begin{theorem}\label{sjogren} For all $n\geq 1$,
$$
F\cap ({\bf r}(n-1)+\mathfrak f^{n+1})=R(n-1)\gamma_{n+1}(F),\ n\geq 1.
$$\end{theorem}

\Para The initial terms  of $\overline E(G)$ can be described by
standard simplicial arguments using the K\"{u}nneth formula and
the Eilenberg-Zilber equivalence. In particular, there exists the
following exact sequence
\begin{equation}\label{tps} 0 \to (H_1(G)\otimes H_2(G))^{\oplus 2}\to \overline
E_{2,\,1}^1(G)\to Tor(H_1(G),\,H_1(G))\to 0
\end{equation}
and, in general,
\begin{equation}\label{abouttor}
0\to T_n(G)\to \overline{E}_{n,\,1}^1(G)\to
Tor_1(\underbrace{H_1(G),\,\dots\,,\, H_1(G)}_{n\ \text{terms}}
)\to 0.
\end{equation}
\par\noindent where $T_1(G)=H_2(G)$ and  $T_{k+1}(G)=H_1(G)\otimes T_k(G)\oplus
T_k(G)\otimes H_1(G),\ k\geq 1$ and for Abelian groups $B_1,\dots,
B_n$, the group $Tor_i(B_1,\dots, B_n)$ denotes the $i$-th
homology group of the complex $P_1\otimes \dots \otimes P_n,$
where $P_j$ is a $\mathbb Z$-flat resolution of $B_j$ for
$j=1,\dots, n$. We clearly have $$Tor_0(B_1,\dots, B_n)=B_0\otimes
\dots \otimes B_n,\ Tor_i(B_1,\dots, B_n)=0,\ i\geq n.$$

\para
For a  description of the initial terms of the spectral sequence
$E(G)$, one needs more complicated theory of derived functors of
polynomial functors. In the quadratic case such a theory was
developed by Baues and Pirashvili \cite{BP:00}; their theory
implies that the terms $E_{2,\,m}^1\ (m\geq 0)$ can be described
explicitly. Let $X$ be a simplicial group which is free Abelian in
each degree. Then there exists [\cite{BP:00}, (4.1)] a natural
short exact sequence of graded Abelian groups
\begin{equation}\label{deri}
0\to Sq^\otimes(\pi_*(X))\rightarrow\pi_*(\wedge^2 X)\rightarrow
Sq^\star(\pi_*(X))[-1] \to 0
\end{equation}
where $\pi_*(X)$ and $\pi_*(\wedge^2 X)$ are the graded homotopy
groups of $X$ and $\wedge^2 X$ respectively. The sequence
(\ref{deri}) gives the following functorial description of the
term $E_{2,\,1}^1$:
\para
 There exists  a natural short exact
sequence:
\begin{equation}\label{lps}
0\to H_1(G)\otimes H_2(G)\to E_{2,\,1}^1(G)\to \Omega(H_1(G))\to
0.
\end{equation}
\para
\noindent{\bf Dimension quotients.} Clearly, for every $n\geq
3$ and $1\leq k\leq n-1$ one has the following commutative diagram
$$
\begin{CD}
E_{n-k,\,1}^k(G) @>d_{n-k,1}^k>> E_{n,\,0}^k(G)
@>>> E_{n,\,0}^{k+1}(G) @>>> 0\\
@V{\kappa_{n-k,\,1}^k}VV @V{\kappa_{n,\,0}^k}VV @V{\kappa_{n,\,0}^{k+1}}VV @.\\
\overline E_{n-k,\,1}^k(G) @>\overline d_{n-k,\,1}^k>> \overline
E_{n,\,0}^k(G) @>>> \overline E_{n,\,0}^{k+1}(G) @>>> 0.
\end{CD}
$$
which has the form
\begin{equation}\label{diahigh1}
\begin{CD}
E_{n-1,\,1}^{1}(G) @>d_{n-1,\,1}^{1}>> E_{n,\,0}^{1}(G)
@>>>  E_{n,\,0}^2(G) @>>> 0\\
@V{\kappa_{n-1,\,1}^{1}}VV @V{\kappa_{n,\,0}^{1}}VV @V{\kappa_{n,\,0}^{2}}VV @.\\
\overline E_{n-1,\,1}^{1}(G) @>\overline d_{n-1,\,1}^{1}>>
\overline E_{n,\,0}^{1}(G) @>>> \overline E_{n,\,0}^2(G) @>>> 0
\end{CD}
\end{equation}
for $k=1$ and
\begin{equation}\label{diahigh}
\begin{CD}
E_{1,\,1}^{n-1}(G) @>d_{1,\,1}^{n-1}>> E_{n,\,0}^{n-1}(G)
@>>> \gamma_n(G)/\gamma_{n+1}(G) @>>> 0\\
@V{\kappa_{1,\,1}^{n-1}}VV @V{\kappa_{n,\,0}^{n-1}}VV @V{\kappa_{n,\,0}^{n}}VV @.\\
\overline E_{1,\,1}^{n-1}(G) @>\overline d_{1,\,1}^{n-1}>>
\overline E_{n,\,0}^{n-1}(G) @>>> \mathfrak g^n/\mathfrak g^{n+1}
@>>> 0
\end{CD}
\end{equation}
for $k=n-1$. Diagrams (\ref{diahigh1}) and (\ref{diahigh})
together with the snake lemma imply the following exact sequences
of Abelian groups:
\begin{equation}\label{seq37}
0 \to \text{Ker}(\kappa_{n,0}^2)\to \text{Coker}(\text{Im}(d_{n-1,1}^1)\to \text{Im}(\overline
d_{n-1,1}^1))\to \text{Coker}(\kappa_{n,0}^1)\to \text{Coker}(\kappa_{n,0}^2)
\to 0 \end{equation} and
\begin{multline}\label{seq38}
0\to \text{Ker}\{\text{Im}(d_{1,1}^{n-1})\to \text{Im}(\overline d_{1,1}^{n-1})\}\to
\text{Ker}(\kappa_{n,0}^{n-1})\to (\gamma_n(G)\cap
D_{n+1}(G))/\gamma_{n+1}(G)\to\\
\text{Coker}\{\text{Im}(d_{1,1}^{n-1})\to \text{Im}(\overline d_{1,1}^{n-1})\}\to
\text{Coker}(\kappa_{n,0}^{n-1})\to \text{Coker}(\kappa_{n,0}^n)\to 0
\end{multline}

As a result we obtain the following diagram which we will use
later: \begin{equation}\label{nachd} \xyma{\text{Ker}(\kappa_{n,\,0}^2)
\ar@{->}[r]  \ar@{>->}[d] & \text{Ker}(\kappa_{n,\,0}^3) \ar@{->}[r] & \dots \ar@{->}[r] & \text{Ker}(\kappa_{n,\,0}^{n-1}) \ar@{->}[d]\\
\text{Coker}(\text{Im}(d_{n-1,1}^1)\to \text{Im}(\overline d_{n-1,1}^1)) & & & \frac{\gamma_n(G)\cap D_{n+1}(G)}{\gamma_{n+1}(G)}\\
}\end{equation}

\para
For an Abelian group $A$, let $L(A)$ denote the free Lie algebra on $A$, and let $L_n(A),\ n\geq 0,$ be its $n$-th homogeneous component; we can view $L_n$ as an endo-functor on the category {\sf Ab} of Abelian groups.  For any endo-functor $F: {\sf Ab} \to {\sf Ab}$, let $\mathcal D_iF,\ i\geq 0,$  denote the derived functor $\mathcal L_iF(-,\,0)$ in the sense of Dold-Puppe \cite{DP:61}.  
Clearly, the map $\kappa_{n-1,1}^1$ can be presented in a natural
diagram
$$
\begin{CD}
0 @>>> K_{n-1}(G) @>>> E_{n-1,\,1}^1(G) @>>> \mathcal D_1L_n(G_{ab}) @>>> 0\\
@. @VVV @V{\kappa_{n-1,\,1}^1}VV @VVV @.\\
0 @>>> T_{n-1}(G) @>>> \overline E_{n-1,\,1}^1(G) @>>>
\mathcal D_1\otimes^n(G_{ab}) @>>> 0
\end{CD}
$$
for some functors $K_{n-1},\ T_{n-1}: {\sf Gr\to Gr}$. We set 
$$^1E_{n,0}^1(G)=\text{Coker}(K_{n-1}(G)\buildrel{d_{n-1,1}^1}\over\to E_{n,0}^1)(G),$$and $$  ^1\overline
E_{n,0}^1(G)=\text{Coker}(T_{n-1}(G)\buildrel{\overline d_{n-1,1}^1}\over\to
\overline E_{n,0}^1(G)).$$ In this notation, we obtain the following
natural diagram:
$$
\begin{CD}
\mathcal D_1L_{n-1}(G_{ab}) @>{^1d_{n-1,1}^1}>> ^1E_{n,0}^1(G) @>>> E_{n,0}^2(G) @>>> 0\\
@VVV @V{^1\kappa_{n,0}^1}VV @V{\kappa_{n,0}^2}VV @.\\
\mathcal D_1\otimes^{n-1}(G_{ab}) @>{^1\overline d_{n-1,1}^1}>> ^1\overline
E_{n,0}^1(G) @>>> \overline E_{n,0}^2(G) @>>> 0
\end{CD}
$$
where the maps $^1d_{n-1,1}^1$ and $^1\overline d_{n-1,1}^1$ are
induced by the maps $d_{n-1,1}^1$ and $\overline d_{n-1,1}^1$
respectively.
\para
Define the functor $S_n: {\sf Ab\to Ab}$ by setting
$$
S_n(A)=\text{Coker}(L_n(A)\to \otimes^n(A)),\ A\in \sf Ab.
$$
Clearly, $\mathcal D_1S_n(A)=\text{Coker}(\mathcal D_1L_n(A)\to \mathcal D_1\otimes^n(A)).$ In this
notation, snake lemma implies the following analog of the diagram
(\ref{nachd}):
\begin{equation}\label{nachd1} \xyma{\text{Ker}( ^1\kappa_{n,\,0}^1) \ar@{->}[d]\\ \text{Ker}(\kappa_{n,\,0}^2)
\ar@{->}[r]  \ar@{->}[dd] & \text{Ker}(\kappa_{n,\,0}^3) \ar@{->}[r] & \dots \ar@{->}[r] & \text{Ker}(\kappa_{n,\,0}^{n-1}) \ar@{->}[d]\\
& & V_{n-1}(G) \ar@{>->}[ld] \ar@{->}[d] \ar@{-->}[llu] &
\frac{\gamma_n(G)\cap D_{n+1}(G)}{\gamma_{n+1}(G)} \ar@{->>}[ld]\\
\text{Coker}(\xi_{n-1}) \ar@{->}[d] & \mathcal D_1S_{n-1}(G_{ab}) \ar@{->>}[l] \ar@{->}[ld] & \frac{\gamma_n(G)\cap D_{n+1}(G)}{\gamma_{n+1}(G).\text{Im}(\text{Ker}( ^1\kappa_{n,\,0}^1))}\\
\text{Coker}( ^1\kappa_{n,0}^1)}\end{equation} where $\xi_{n-1}:
\text{Im}(^1d_{n-1,1}^1)\to \text{Im}(^1\overline d_{n-1,1}^1)$ and $V_{n-1}(G)$
is the kernel of the composition of the natural maps
$\mathcal D_1S_{n-1}(G_{ab})\to \text{Coker}(\xi)\to \text{Coker}(^1\kappa_{n,\,0}^1)$. 
[The dotted  arrow from $A$ to $B$ in the diagram above, and also later in diagram \ref{bigdiagram},  means that there is a map $A\to B/C$ for some subgroup C.]
Hence, for every group $G$ and $n\geq 3$, there is a natural
subgroup of $\mathcal D_1S_{n-1}(G_{ab})$, namely $V_{n-1}(G)$, which maps
canonically to the quotient $\frac{\gamma_n(G)\cap
D_{n+1}(G)}{\gamma_{n+1}(G).\text{Im}(\text{Ker}(^1\kappa_{n,\,0}))}.$ Denote
this map by
$$
v_n: V_{n-1}(G)\to \frac{\gamma_{n}(G)\cap
D_{n+1}(G)}{\gamma_{n+1}(G).\text{Im}(\text{Ker}(^1\kappa_{n,\,0}^1))}
$$

\par\vspace{.25cm}\noindent {\bf Remark.} It may be noted that, heuristically speaking,
the contribution of the functor $\mathcal D_1S_n(G_{ab})$ to the structure
of the $(n+1)$-st dimension quotient becomes smaller and smaller
as $n$ increases.

\Para {\bf Remark.} The map $\overline d_{n-1,1}^1: T_{n-1}(G)\to
\overline E_{n,0}^1(G)$ can be described explicitly (see, for
example, the description of $\overline d_{1,1}^1$ in the Section
\ref{sew}) and the exact sequence (\ref{abouttor}) implies that
$$
\text{Im}(\text{Ker}(^1\kappa_{n,\,0}^1))=R\gamma_n(F)\cap (1+(R\cap
F'-1)(n-2)+\mathfrak f^{n+1})/R\gamma_{n+1}(F)
$$
In particular, $\text{Im}(\text{Ker}(^1\kappa_{n,\,0}^1))=0$, provided $H_2(G)=0$
by Theorem \ref{sjogren}.
\\

We analyse the map $v_n$ for the case $n=3$. Observe that
$$\mathcal D_1S_2(G_{ab})= G_{ab} \stackrel{\wedge}{*}G_{ab}.$$ We need an
identification theorem which we present in the next section.
\Para
\section{An identification theorem}\label{it}
\para
\begin{theorem}\label{msq} If  $F$ is a free group and $R$ a normal subgroup of $F$, then
$$F\cap(1+\mathfrak f(R\cap F'-1)+(R\cap F'-1)\mathfrak f+{\bf
r}(2)+{\mathfrak f}^4)=[R\cap F',\,F]\gamma_4(F), $$ where $F'$ is
the derived subgroup of $F$.
\end{theorem}\para
We need the following
\begin{lemma}\label{idlemma}
Let $F$ be a free group with a basis $\{x_1,\,\dots,\, x_m\}$, $u$ an
element of \linebreak  $(F'-1)\mathfrak f,$ such that
$$
w-1\equiv u+c.v\mod {\mathfrak f}^4,
$$
for some $c>0$ and $w\in \gamma_3(F),\ v\in {\mathfrak f}^3$. Then
$$
u\equiv c.v_1\mod {\mathfrak f}^4,
$$
where $v_1\in (F'-1)\mathfrak f.$
\end{lemma}
\begin{proof}
We can view ${\mathfrak f}^3/{\mathfrak f}^4$ as a free Abelian
group isomorphic to $F_{ab}^{\otimes 3}$ with a basis
$\{x_i\otimes x_j\otimes x_k\ |\ i,\,j,\,k=1,\,\dots,\,m\}.$
Modulo ${\mathfrak f}^4,$ the group $(F'-1)\mathfrak f$ is
generated by elements
$$
([x_i,\,x_j]-1)(x_k-1),\ i,\,j,\,k=1,\,\dots,\,m.
$$
Let
$$
u\equiv\sum_{i,\,j,\,k}d_{i,\,j,\,k}([x_i,\,x_j]-1)(x_k-1)\mod
{\mathfrak f}^4.
$$
For a given triple $(i,\,j,\,k)$, the sum of all coefficients in
$u$, which contribute to $x_i\otimes x_j\otimes x_k$ (and
$S_3$-permutations) must be divided by $c$. We have in ${\mathfrak
f}^3/{\mathfrak f}^4:$
\begin{align*}
& ([x_i,\,x_j]-1)(x_k-1)\mapsto x_i\otimes x_j\otimes x_k-x_j\otimes
x_i\otimes x_k,\\
& ([x_k,\,x_i]-1)(x_j-1)\mapsto x_k\otimes x_i\otimes x_j-x_i\otimes
x_k\otimes x_j,\\
& ([x_j,\,x_k]-1)(x_i-1)\mapsto x_j\otimes x_k\otimes x_i-x_k\otimes
x_j\otimes x_i
\end{align*}
and no more terms can contribute to $x_i\otimes x_j\otimes x_k$
(and permutations). Clearly  only the product of commutators  of
the form $[x_i,\,x_j,\,x_k]$ and $[x_k,\,x_i,\,x_j]$ with different powers
can contribute from the element $w$. We have:
\begin{align*}
& 1-[x_i,\,x_j,\,x_k]^{f_{i,\,j,\,k}}\mapsto f_{i,\,j,\,k}(x_i\otimes
x_j\otimes x_k-x_j\otimes x_i\otimes x_k+x_k\otimes x_j\otimes
x_i-x_k\otimes x_i\otimes x_j)\\
& 1-[x_k,\,x_i,\,x_j]^{f_{k,\,i,\,j}}\mapsto f_{k,\,i,\,j}(x_k\otimes
x_i\otimes x_j-x_i\otimes x_k\otimes x_j+x_j\otimes x_i\otimes
x_k-x_j\otimes x_k\otimes x_i)
\end{align*}
Now we have
\begin{multline}\label{gjk}
d_{i,\,j,\,k}(x_i\otimes x_j\otimes x_k-x_j\otimes x_i\otimes
x_k)+\\ d_{k,\,i,\,j}(x_k\otimes x_i\otimes x_j-x_i\otimes x_k\otimes
x_j)+  d_{j,\,k,\,i}(x_j\otimes x_k\otimes x_i-x_k\otimes x_j\otimes
x_i)+\\
f_{i,\,j,\,k}(x_i\otimes x_j\otimes x_k-x_j\otimes x_i\otimes
x_k+x_k\otimes x_j\otimes x_i-x_k\otimes x_i\otimes
x_j)+\\f_{k,\,i,\,j}(x_k\otimes x_i\otimes x_j-x_i\otimes x_k\otimes
x_j+x_j\otimes x_i\otimes x_k-x_j\otimes x_k\otimes x_i)\\
= x_i\otimes x_j\otimes x_k(d_{i,\,j,\,k}+f_{i,\,j,\,k})+x_j\otimes x_i\otimes x_k(-d_{i,\,j,\,k}-f_{i,\,j,\,k}+f_{k,\,i,\,j})+\\
x_k\otimes x_i\otimes
x_j(d_{k,\,i,\,j}-f_{i,\,j,\,k}+f_{k,\,i,\,j})+x_i\otimes x_k\otimes
x_j(-d_{k,\,i,\,j}-f_{k,\,i,\,j})+\\
x_j\otimes x_k\otimes x_i(d_{j,\,k,\,i}-f_{k,\,i,\,j})+x_k\otimes
x_j\otimes x_i(-d_{j,\,k,\,i}+f_{i,\,j,\,k})=\\ c.v(x_i,\,x_j,\,x_k),
\end{multline}
for some element $v(x_i,\,x_j,\,x_k)$ from ${\mathfrak
f}^3/{\mathfrak f}^4$. Suppose first that all $i,\,j,\,k$ are
different. Then (\ref{gjk}) implies that
\begin{align*} d_{i,\,j,\,k}+f_{i,\,j,\,k},\
-d_{i,\,j,\,k}-f_{i,\,j,\,k}+f_{k,\,i,\,j},\
d_{k,\,i,\,j}-f_{i,\,j,\,k}+f_{k,\,i,\,j},\\
-d_{k,\,i,\,j}-f_{k,\,i,\,j},\ d_{j,\,k,\,i}-f_{k,\,i,\,j},\
-d_{j,\,k,\,i}+f_{i,\,j,\,k}\end{align*} must be divided by $c$.
Therefore, all $d_{i,\,j,\,k}, d_{k,\,i,\,j},\, d_{j,\,k,\,i},\,
f_{i,\,j,\,k}, f_{k,\,i,\,j}$ must be divisible  by $c$.
\para
Suppose $i=j$. Then we reduce all the expression to the element
$d_{i,\,k,\,i}([x_i,\,x_k]-1)(x_i-1)$ and a bracket
$1-[x_i,\,x_k,\,x_i]^{f_{i,\,k,\,i}}$. We then have
$$
-f_{i,\,k,\,i}x_i\otimes x_i\otimes
x_k+(2f_{i,\,k,\,i}+d_{i,\,k,\,i})x_i\otimes x_k\otimes
x_i+(-f_{i,\,k,\,i}-d_{i,\,k,\,i})x_k\otimes x_i\otimes x_i=c.v(x_i,\,x_k)
$$
and we again conclude that $d_{i,\,k,\,i}$ and $f_{i,\,k,\,i}$ must be
divisible  by $c$.
\end{proof}
\par\vspace{.5cm}
Next, let $w\in \gamma_3(F)$ and
$$
w-1\equiv u\mod {\mathfrak f}^4,
$$
where
$$
u\in (R\cap F'-1)\mathfrak f+\mathfrak f(R\cap F'-1)+{\mathfrak
r}(2).
$$
We claim that
$$
w=w_1w_2\mod [R\cap F',\,F][R,\,F,\,F]\gamma_4(F),
$$
such that
\begin{align*}
& w_1-1\in (R\cap F'-1)\mathfrak f+\mathfrak f(R\cap
F'-1)+{\mathfrak
f}^4,\\
& w_2-1\in {\bf r}(2)+{\mathfrak f}^4,
\end{align*}
that is, we can divide the problem of identification of the
subgroup
$$
F\cap (1+(R\cap F'-1)\mathfrak f+\mathfrak f(R\cap F'-1)+{\bf
r}(2)+{\mathfrak f}^4)
$$
into two parts:

\para\noindent
$(i)$ identification of
$$
F\cap (1+(R\cap F'-1)\mathfrak f +\mathfrak f(R\cap
F'-1)+{\mathfrak f}^4)
$$
and
\para\noindent
$(ii)$ identification of
$$
F\cap (1+{\bf r}(2)+{\mathfrak f}^4)
$$
Let $u=u_1+u_2,$ where
$$
u_1\in (R\cap F'-1) \mathfrak f+\mathfrak f(R\cap F'-1),\ u_2\in
{\bf r}(2).
$$
Since we can work modulo $[R\cap F',\,F]$, we can assume that
$$
u_1\in (R\cap F'-1)\mathfrak f.
$$
Now, using the argument of Gupta \cite[Lemma 1.5(B), p.\.72]{Gupta:87}, we can easily
conclude that
$$
u_2\equiv e_m.v_2\mod {\mathfrak f}^4,
$$
where $v_2$ involves the element $x_m$. Hence, by Lemma \ref{idlemma},
we conclude that $u_1=e_m.v_3,$ where $v_3$ contains some
nontrivial entries of the element $x_m$. Therefore, we have
$$
u=e_m.v \mod {\mathfrak f}^4,
$$
where $v$ is an element from $(R\cap F'-1)\mathfrak f+\mathfrak
f(R\cap F'-1)+{\bf r}(2).$ Since $u$ is a Lie element, $v$ is
again a Lie element and we conclude that
$$
w\equiv w'^{e_m}\mod [R\cap F',\,F][R,\,F,\,F]\gamma_4(F).
$$
Now we can delete all the brackets from $w'$ with entries of
$x_m$ and make an induction. The induction argument shows the
following: 
\begin{quote}{\it Let $w\in \gamma_3(F)$, such that
$$
w-1\in (R\cap F'-1)\mathfrak f+\mathfrak f(R\cap F'-1)+{\bf
r}(2)+{\mathfrak f}^4,
$$
then
$$
w=w_1w_2\mod [R\cap F',\,F][R,\,F,\,F]\gamma_4(F),
$$
such that
\begin{align*}
& w_1-1\in (R\cap F'-1)\mathfrak f+\mathfrak f(R\cap
F'-1)+{\mathfrak
f}^4,\\
& w_2-1\in {\bf r}(2)+{\mathfrak f}^4
\end{align*}
} \end{quote}
and Theorem \ref{msq} follows.

\Para
\section{The fourth dimension quotient (cont'd.)}\label{sew}

\para
Since $D_4(G)\subseteq \gamma_3(G),$ the dimension quotient
$D_4(G)/\gamma_4(G)$ is exactly the kernel of the map
\begin{equation}
\kappa_{3,\,0}^3: E_{3,\,0}^3(G)\to \overline
E_{3,\,0}^3(G).\end{equation} For $n=3$ the sequence (\ref{seq38})
reduces to the following
$$
0\to \text{Ker}(\text{Im}(d_{1,\,1}^2)\to \text{Im}(\overline d_{1,\,1}^2))\to
\text{Ker}(\kappa_{3,\,0}^2)\to D_4(G)/\gamma_4(G)\to 1
$$
for every group $G$. The sequence (\ref{seq37}) has the following
form:
$$
0\to \text{Ker}(\kappa_{3,\,0}^2)\to \text{Coker}(\eta)\to
\text{Coker}(\kappa_{3,\,0}^1)\to \text{Coker}(\kappa_{3,\,0}^2)\to 0,
$$
where $\eta: \text{Im}(d_{1,1}^1)\to \text{Im}(\overline d_{1,1}^1).$ From the
exact sequences (\ref{tps}) and (\ref{lps}), we have
 a commutative diagram
\begin{equation}
\begin{CD}
0 @>>> H_1(G)\otimes H_2(G) @>>> E_{2,\,1}^1(G) @>>> \Omega(H_1(G))
@>>> 0\\
@. @VVV @V{\kappa_{2,\,1}^1}VV @V{T}VV @.\\
0 @>>> (H_1(G)\otimes H_2(G))^{\oplus 2} @>>> \overline
E_{2,\,1}^1(G) @>>> \text{Tor}(H_1(G),\,H_1(G)) @>>> 0
\end{CD}
\end{equation}
where $T$ is as in (\ref{eilmac1}).

\para

Every element $x\in H_2(G)=H_2(F/R)$ can be presented in the form
$x\equiv \prod_{i=1}^k[f_1^{(i)},\,f_2^{(i)}] \mod \gamma_3(F),$
with $\prod_{i=1}^k[f_1^{(i)},\,f_2^{(i)}]\in R.$ Then the map
$d_{1,\,1}^1$ restricted to the component $H_2(G)\otimes H_1(G)$
is given by
$$
x\otimes \bar g\mapsto \prod_{i=1}^k
[f_1^{(i)},f_2^{(i)},g].[R,F,F]\gamma_4(F),\ x\in H_2(G),\ \bar
g\in G_{ab},
$$
where $\bar f$ is the image of $f\in F$ in $G_{ab}$. On the other
hand, it is easy to see that the map $\overline d_{1,1}^1$
restricted to $(H_2(G)\otimes H_1(G))^{\oplus 2}$ is induced by
$$
(x\otimes \bar g_1,\  x \otimes \bar g_2)\mapsto \sum_{i=1}^k \bar
f_1^{(i)}\otimes \bar f_2^{(i)}\otimes \bar g_1+\sum_{i=1}^k \bar
g_2\otimes \bar f_1^{(i)}\otimes \bar f_2^{(i)},\ x\in H_2(G),\
\bar g_1,\, \bar g_2\in G_{ab}. $$ Hence we have the following
diagrams:
$$
\begin{CD}
0 @>>> H_1(G)\otimes H_2(G)@>>> E_{1,\,1}^1(G) @>>> \Omega(H_1(G))@>>> 0\\
@. @VVV @V{d_{1,\,1}^1}VV @VVV @.\\
0 @>>> \frac{[F,\,R\cap
\gamma_2(F)]\gamma_4(F)}{[R,\,F,\,F]\gamma_4(F)} @>>>
\frac{\gamma_3(F)}{[R,\,F,\,F]\gamma_4(F)} @>>>
\frac{\gamma_3(F)}{[F,\,R\cap \gamma_2(F)]\gamma_4(F)} @>>> 0\\
@. @. @VVV @VVV @.\\
@. @. E_{3,0}^2(G) @>{\simeq}>> E_{3,0}^2(G)
\end{CD}
$$
$$
\begin{CD}
0 @>>> (H_1(G)\otimes H_2(G))^{\oplus 2}@>>> \overline E_{1,\,1}^1(G) @>>> Tor(H_1(G),H_1(G))@>>> 0\\
@. @VVV @V{\overline d_{1,\,1}^1}VV @VVV @.\\
0 @>>> \frac{(R\cap \gamma_2(F)-1)(1)+\mathfrak r(2)+\mathfrak
f^4}{\mathfrak r(2)+\mathfrak f^4} @>>> \frac{\mathfrak
f^3}{\mathfrak r(2)+\mathfrak f^4} @>>>
\frac{\mathfrak f^3}{(R\cap \gamma_2(F)-1)(1)+\mathfrak r(2)+\mathfrak f^4} @>>> 0\\
@. @. @VVV @VVV @.\\
@. @. \overline E_{3,0}^2(G) @>{\simeq}>> \overline E_{3,0}^2(G)
\end{CD}
$$
and
\begin{equation}\label{va1}
\begin{CD}
\Omega(G_{ab}) @>{^1d_{1,\,1}^{1}}>> \frac{\gamma_3(F)}{[R\cap \gamma_2(F),\,F]\gamma_4(F)}@>>> E_{3,\,0}^2(F) @>>> 0\\
@V{T}VV @V{^1\kappa_{3,\,0}^1}VV @V{\kappa_{3,\,0}^2}VV @.\\
Tor(G_{ab},\,G_{ab}) @>{^1\overline d_{1,\,1}^{1}}>>
\frac{\mathfrak f^3}{(R\cap F'-1)(1)+\mathfrak r(2)+\mathfrak
f^4}@>>> \overline E_{3,\,0}^2(F) @>>> 0
\end{CD}
\end{equation}
where $^1d_{1,\,1}^{1}$ and $^1\overline d_{1,\,1}^{1}$ are
induced by $d_{1,\,1}^1$ and $\overline d_{1,\,1}^1$ respectively; here, and subsequently, $(R\cap F'-1)(1):= (R\cap F'-1)\mathfrak f+\mathfrak f(R\cap F'-1)$. 
\para
Now Theorem \ref{msq} implies the following commutative diagram
\begin{equation}\label{poew33}
\xyma{ & \text{Im}(^1d_{1,1}^{1}) \ar@{>->}[r] \ar@{>->}[d]_{\eta'} &
\frac{\gamma_3(F)}{[R\cap \gamma_2(F),\,F]\gamma_4(F)}
\ar@{->>}[r] \ar@{>->}[d]_{^1\kappa_{3,\,0}^1}
& E_{3,\,0}^2(F) \ar@{->}[d]_{\kappa_{3,0}^2}\\
& \text{Im}(^1\overline d_{1,1}^{1}) \ar@{>->}[r] \ar@{->>}[d] &
\frac{\mathfrak f^3}{(R\cap F')(1)+\mathfrak r(2)+\mathfrak f^4}
\ar@{->>}[r] \ar@{->>}[d] & \overline
E_{3,\,0}^2(F)\\
\text{Ker}(\kappa_{3,0}^2) \ar@{>->}[r] & \text{Coker}(\eta') \ar@{->}[r] &
\text{Coker}(^1\kappa_{3,\,0}^1)\\
& G_{ab} \stackrel{\wedge}{*}G_{ab} \ar@{->>}[u] \ar@{->}[ur]\\
V_2(G) \ar@{->>}[uu] \ar@{>->}[ur]
 }
\end{equation}
where $\eta'$ is the map induced by $T$.
\para
Hence for $n=3$ the diagram \ref{nachd} has a simple form and we have the following:

\begin{theorem} There exists the following natural system of monomorphism and epimorphisms
$$
\xyma{\text{Ker}(\kappa_{3,\,0}^2)
\ar@{->>}[r] & D_4(G)/\gamma_4(G)\\
V_2(G) \ar@{->>}[u] \ar@{->>}[ru] \ar@{>->}[r] & G_{ab} \stackrel{\wedge}{*}G_{ab}  \\
}
$$
\end{theorem}

As an immediate consequence of the above result, we have various  conditions for the fourth dimension quotient of a group $G$ to be trivial.
\para
\begin{cor}\label{4dim}
If either $V_2(G),\ or \  G_{ab} \stackrel{\wedge}{*}G_{ab} $ is
trivial, then $D_4(G)=\gamma_4(G).$
\end{cor}
\par\vspace{.25cm}\noindent{\bf Remark.} It may be noted that triviality of the exterior square $ G_{ab} \stackrel{\wedge}{*}G_{ab}$ , for finitley generated $G$, means that all primary parts of $G_{ab}$ are cyclic.

\Para
\section{The fifth dimension quotient}\label{fifthdq}

\Para For $n=4$, the diagram (\ref{nachd1}) leads to the following:
\begin{equation}\label{fifth1}
\xyma{& \text{Ker}(^1\kappa_{4,\,0}^1) \ar@{->}[d]\\
& \text{Ker}(\kappa_{4,\,0}^2)
\ar@{->}[rr]  \ar@{->}[d] & & \text{Ker}(\kappa_{4,\,0}^3) \ar@{->}[d] \ar@{->}[r] & \text{Coker}(\eta)\\
D_1S_{3}(G_{ab}) \ar@{->>}[r] & \text{Coker}(\xi_4)  & & \frac{\gamma_4(G)\cap D_{5}(G)}{\gamma_{5}(G)} \ar@{->}[d]\\
& & & \text{Coker}\{\text{Im}(d_{1,1}^3)\to \text{Im}(\overline d_{1,1}^3)\} \\
}\end{equation} Where $\eta: \text{Im}(d_{2,1}^2)\to \text{Im}(\overline
d_{2,1}^2)$. \Para Let us examine the other derived functors for
their contribution to the dimension quotients. We have the
following diagram:
$$
\xyma{ \text{Ker}(^1d_{1,1}^1) \ar@{->}[r] \ar@{>->}[d] & \Omega(G_{ab})
\ar@{->}[r]^{^1d_{1,1}^1} \ar@{>->}[d] & ^1E_{3,0}^1 (G)\ar@{->}[d]\\
\text{Ker}(^1\overline d_{1,1}^1) \ar@{->}[r] & Tor(G_{ab},\,G_{ab})
\ar@{->}[r]^{^1\overline d_{1,1}^1} & ^1\overline E_{3,0}^1(G)\\
}
$$
and natural quotients
$$
\xyma{ H_3(G) \ar@{->}[r] \ar@{=}[d] & \text{Ker}(^1d_{1,1}^1)
\ar@{->>}[r]
\ar@{>->}[d] & \overline{\text{Ker}(^1d_{1,1}^1)} \ar@{>->}[d]\\
H_3(G) \ar@{->}[r] & \text{Ker}(^1\overline d_{1,1}^1) \ar@{->>}[r]
\ar@{->>}[d] &
\overline{\text{Ker}(^1\overline d_{1,1}^1)} \ar@{->>}[d]\\
& \text{Ker}(^1\overline d_{1,1}^1)/\text{Ker}(^1d_{1,1}^1) \ar@{=}[r] &
\overline{\text{Ker}(^1\overline d_{1,1}^1)}/\overline{\text{Ker}(^1d_{1,1}^1)}\\
}
$$
Now we have the following commutative diagram:
$$
\xyma{ \overline{\text{Ker}(^1d_{1,1}^1)} \ar@{->}[r]^{d_{2,1}^2} \ar@{>->}[d] & E_{4,0}^2(G) \ar@{->}[d]_{\kappa_{4,0}^2}\\
\overline{\text{Ker}(^1\overline d_{1,1}^1)} \ar@{->}[r]^{\overline
d_{2,1}^2} & \overline
E_{4,0}^2(G)\\
}
$$
Analogous to the diagram (\ref{va1}), we have  the following
commutative diagram
\begin{equation}\label{va2}
\begin{CD}
\overline{\text{Ker}(^1d_{1,\,1}^1)} @>{^1d_{2,\,1}^{2}}>> \frac{\gamma_4(F)}{([R\cap \gamma_2(F),\,F]\cap \gamma_4(F))\gamma_5(F)}@>>> E_{4,\,0}^3(G) @>>> 0\\
@V{T''}VV @V{^1\kappa_{4,0}^2}VV @V{\kappa_{4,\,0}^2}VV @.\\
\overline{\text{Ker}(^1\overline d_{1,\,1}^1)} @>{^1\overline
d_{2,\,1}^{2}}>> \frac{\mathfrak f^4}{(R\cap F'-1)(1)\cap
\mathfrak f^4+\mathfrak r(2)\cap\mathfrak f^4+\mathfrak f^5}@>>>
\overline E_{4,\,0}^3(G) @>>> 0
\end{CD}
\end{equation}
where the maps $^1d_{2,\,1}^2,$ $^1\overline d_{2,\,1}^2,$
$^1\kappa_{4,0}^2$ and $T''$ are induced by the maps
$d_{2,\,1}^2,$ $\overline d_{2,\,1}^2,$ $\kappa_{4,0}^2$ and $T$
respectively. Let
$$ V(F,R):=F\cap (1+(R\cap F'-1)(1)\cap \mathfrak f^4+{\bf
r}(2)\cap\mathfrak f^4+\mathfrak f^5).
$$
An application of the snake lemma  implies the following exact sequence:
$$
\frac{V(F,R)}{([R\cap \gamma_2(F),\,F]\cap
\gamma_4(F))\gamma_5(F)}\to \text{Ker}(\kappa_{4,0}^3)\to
\text{Coker}(\text{Im}(^1d_{2,1}^{2})\to \text{Im}(^1\overline d_{2,1}^{2}))
$$
with a natural sequence of epimorphisms and monomorphisms:
\begin{equation}\label{fifth2}
\xyma{\text{Ker}(^1d_{1,1}^1)/\text{Ker}(^1\overline d_{1,1}^1)
\ar@{>->}[r] \ar@{=}[d] &  G_{ab} \stackrel{\wedge}{*}G_{ab}\\
\overline{\text{Ker}(^1d_{1,1}^1)}/\overline{\text{Ker}(^1\overline d_{1,1}^1)} \ar@{->>}[r] & \text{Coker}(\text{Im}(^1d_{2,1}^{2})\to \text{Im}(^1\overline d_{2,1}^{2}))\\
}\end{equation}

We collect the above diagrams into the following:
\newpage
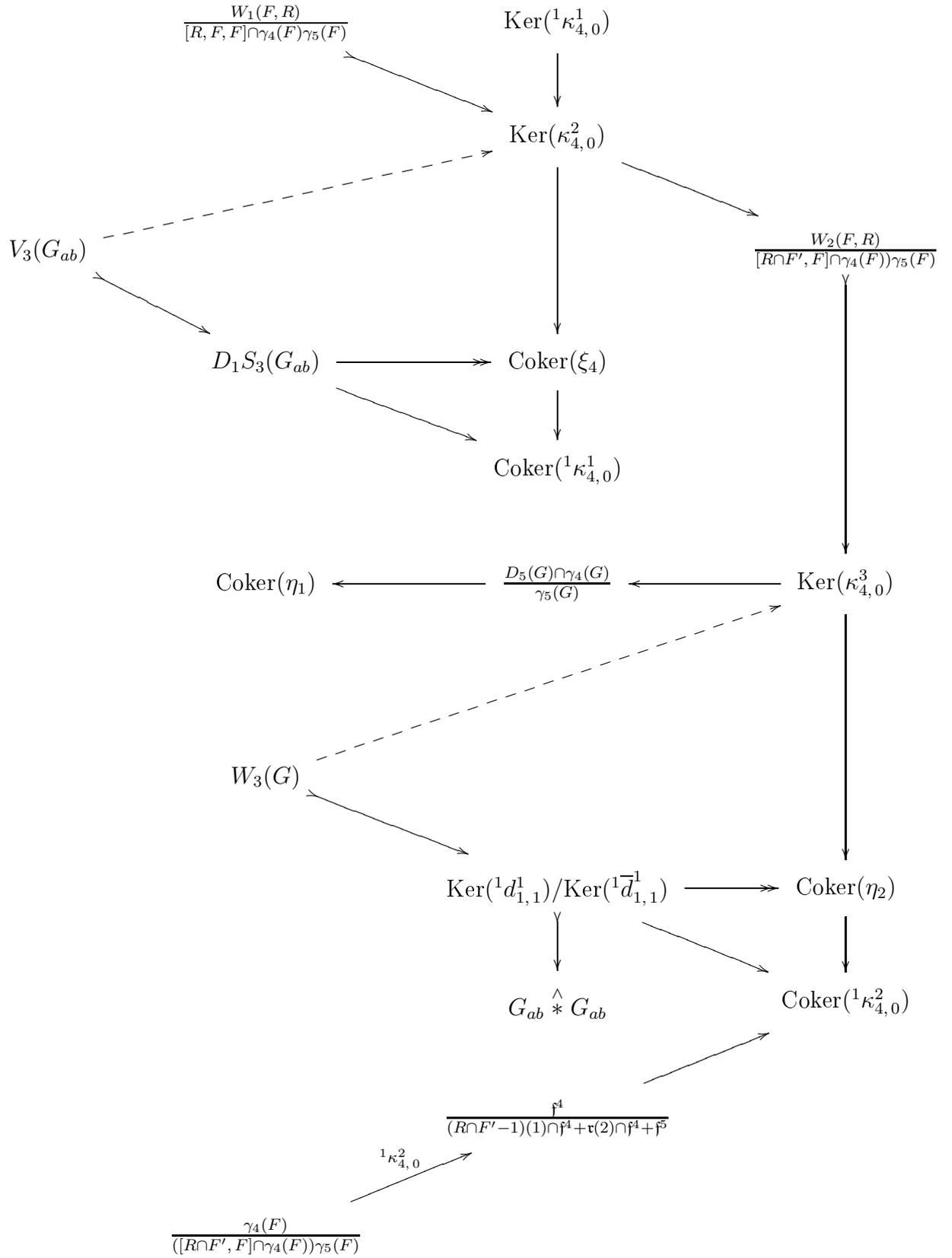
\begin{equation}\label{bigdiagram}
\xyma{& \frac{W_1(F,\,R)}{[R,\,F,\,F]\cap \gamma_4(F)\gamma_5(F)} \ar@{>->}[dr] & \text{Ker}(^1\kappa_{4,\,0}^1) \ar@{->}[d]\\
& & \text{Ker}(\kappa_{4,\,0}^2)
\ar@{->}[dd] \ar@{->}[dr]\\
  V_3(G_{ab}) \ar@{>->}[dr] \ar@{-->}[urr] & & & \frac{W_2(F,\,R)}{[R\cap F',\,F]\cap \gamma_4(F))\gamma_5(F)} \ar@{>->}[ddd]\\
 & D_1S_3(G_{ab}) \ar@{->>}[r] \ar@{->}[rd] & \text{Coker}(\xi_4) \ar@{->}[d]\\
 & & \text{Coker}(^1\kappa_{4,\,0}^1)\\ & \text{Coker}(\eta_1) & \frac{D_5(G)\cap \gamma_4(G)}{\gamma_5(G)} \ar@{->}[l] & \text{Ker}(\kappa_{4,\,0}^3) \ar@{->}[l] \ar@{->}[ddd]\\
\\
& W_3(G) \ar@{>->}[rd] \ar@{-->}[rruu] \\
& & \text{Ker}(^1d_{1,\,1}^1)/\text{Ker}(^1\overline d_{1,\,1}^1) \ar@{->>}[r] \ar@{->}[rd] \ar@{>->}[d] & \text{Coker}(\eta_2)\ar@{->}[d]\\
& & G_{ab} \stackrel{\wedge}{*}G_{ab} & \text{Coker}(^1\kappa_{4,\,0}^2)\\
& & \frac{\mathfrak f^4}{(R\cap F'-1)(1)\cap \mathfrak
f^4+\mathfrak r(2)\cap \mathfrak f^4+\mathfrak f^5}\ar@{->}[ur]\\
& \frac{\gamma_4(F)}{([R\cap
F',\,F]\cap\gamma_4(F))\gamma_5(F)}\ar@{->}[ru]^{^1\kappa_{4,\,0}^2}
}\end{equation}

Here \begin{align*}& W_1(F,\,R)=[R\cap F',\,F]\cap \gamma_4(F)\cap
(1+{\bf r}(2)\cap \mathfrak f^4+\mathfrak f^5)\gamma_5(F)\\ &
W_2(F,\,R)=F\cap (1+(R\cap F'-1)(1)\cap \mathfrak f^4+{\bf
r}(2)\cap\mathfrak f^4+\mathfrak f^5)\\ & \eta_1: \text{Im}(d_{1,\,1}^3)\to
\text{Im}(\overline d_{1,\,1}^3),\\
& \eta_2: \text{Im}(^1d_{2,\,1}^2)\to \text{Im}(^1\overline d_{2,\,1}^2)
\end{align*}
and $W_3(G)$ is the kernel of the composition of the natural maps
$$
\text{Ker}(^1d_{1,\,1}^1)/\text{Ker}(^1\overline d_{1,\,1}^1)\to
\text{Coker}(\eta_2)\to \text{Coker}(^1\kappa_{4,\,0}^2).
$$

\para
Suppose now that $H_2(G)=0$. Then $R\cap F'=[R,\,F],$ hence
\begin{align*}
& (R\cap F'-1)(1)\subseteq {\bf r}(2),\\ & W_1(F,\,R)=[F,\,R,\,R]\cap
\gamma_4(F)\gamma_5(F)
\end{align*}
and
$$
\text{Ker}(^1\kappa_{n,\,0}^1)=0
$$
for all $n\geq 3$. Also $\eta_1$ is an isomorphism. Therefore,
diagram (\ref{bigdiagram}) implies the following diagram with exact
horizontal sequence:
$$
\xyma{& V_3(G) \ar@{->>}[d] & & G_{ab} \stackrel{\wedge}{*}G_{ab} & \text{Ker}(^1d_{1,\,1}^1)/\text{Ker}(^1\overline d_{1,\,1}^1) \ar@{>->}[l] \ar@{->>}[d]\\
& \text{Ker}(\kappa_{4,\,0}^2)
\ar@{->}[rr]  & & \text{Ker}(\kappa_{4,\,0}^3) \ar@{->>}[d] \ar@{->}[r] & \text{Coker}(\eta_4)\\
& &  & \frac{\gamma_4(G)\cap D_{5}(G)}{\gamma_{5}(G)}\\
}$$ Thus, in particular, we obtain the following\Para
\begin{theorem}\label{5dim}
Let $G$ be a group with $H_2(G)=0,\ G_{ab}
\stackrel{\wedge}{*}G_{ab}=0,\ \mathcal D_1S_3(G_{ab})=0$. Then
$D_5(G)=\gamma_5(G)$.
\end{theorem}

\bibliographystyle{alpha}


\Para
M. Hartl\\
LAMAV, FR CNRS 2956, ISTV2\\
Universit\'{e} de Valenciennes et du Hainaut-Cambrr\'{e}sis\\Le Mont Houy, 59313 Valenciennes, Cedex 9, France\\
email: Manfred.Hartl@univ-valenciennes.fr
\par\vspace{1cm}\noindent
R. Mikhailov\\
Steklov Mathematical Institute\\
Department of Algebra\\
Gubkina 8\\
Moscow 119991, Russia\\
email: romanvm@mi.ras.ru
\Para
I. B. S. Passi\\
Centre for Advanced Study in Mathematics\\
Panjab University
\\
Chandigarh 160014,
India\\
and \\
Indian Institute of Science Education and Research, Mohali\\
Transit campus: MGSIPA Complex, Sector 26\\
Chandigarh 160019, India\\
email: ibspassi@yahoo.co.in
\end{document}